\documentclass[12pt,a4paper,twoside]{article}

\usepackage{hyperref}

\newcommand{\pageformat}[6]{\setlength{\hoffset}{-1in}
                  \setlength{\voffset}{-1in}
                  \addtolength{\hoffset}{#5}
                            \addtolength{\voffset}{#6}
                            \setlength{\oddsidemargin}{#1}
                            \setlength{\evensidemargin}{#2}
                            \setlength{\textwidth}{\paperwidth}
                  \addtolength{\textwidth}{-\oddsidemargin}
                  \addtolength{\textwidth}{-\evensidemargin}
                  \addtolength{\textwidth}{-\marginparsep}
                  \addtolength{\textwidth}{-\marginparwidth}
                            \setlength{\topmargin}{#3}
                            \setlength{\textheight}{\paperheight}
                  \addtolength{\textheight}{-\topmargin}
                  \addtolength{\textheight}{-\headheight}
                  \addtolength{\textheight}{-\headsep}
                  \addtolength{\textheight}{-\footskip}
                  \addtolength{\textheight}{-#4}}
\pageformat{2cm}{3cm}{25mm}{25mm}{1pt}{0pt}

\usepackage{ifthen}
\newboolean{article}
    \setboolean{article}{true}
\newboolean{report}
\newboolean{book}
\newboolean{letter}
\newboolean{german}
\newboolean{italian}
\newboolean{nobaselinestretch}
\newboolean{nosectionappendix}
\newboolean{oldtoc}
\newboolean{nosectionequn}
\newboolean{notheorem}

\ifthenelse{\boolean{german}}{
    \usepackage{german}}{}

\usepackage[latin1]{inputenc}

\usepackage{amsmath}
\usepackage{amssymb}
\usepackage[mathscr]{eucal}

\ifthenelse{\boolean{notheorem}}{}{
    \usepackage{theorem}}

\ifthenelse{\boolean{nobaselinestretch}}{}{
    \renewcommand{\baselinestretch}{1.25}}

\newenvironment{env}[2]{\begin{#1}#2\end{#1}}{}
    \newcommand{\beq}[1]{\begin{env}{equation}{#1}}
    \newcommand{\beqn}[1]{\begin{env}{equation*}{#1}}
    \newcommand{\bal}[1]{\begin{env}{align}{#1}}
    \newcommand{\baln}[1]{\begin{env}{align*}{#1}}
    \newcommand{\bga}[1]{\begin{env}{gather}{#1}}
    \newcommand{\bgan}[1]{\begin{env}{gather*}{#1}}
    \newcommand{\bflal}[1]{\begin{env}{flalign}{#1}}
    \newcommand{\bflaln}[1]{\begin{env}{flalign*}{#1}}
    \newcommand{\bmu}[1]{\begin{env}{multline}{#1}}
    \newcommand{\bmun}[1]{\begin{env}{multline*}{#1}}
    \newcommand{\bsp}[1]{\begin{env}{split}{#1}}

    \newcommand{\eeq}{\end{env}}
    \newcommand{\eeqn}{\end{env}}
    \newcommand{\eal}{\end{env}}
    \newcommand{\ealn}{\end{env}}
    \newcommand{\ega}{\end{env}}
    \newcommand{\egan}{\end{env}}
    \newcommand{\eflal}{\end{env}}
    \newcommand{\eflaln}{\end{env}}
    \newcommand{\emu}{\end{env}}
    \newcommand{\emun}{\end{env}}
    \newcommand{\esp}{\end{env}}

\newcommand{\lf}{\vspace{2ex}}
\newcommand{\bulletline}[1][]{\lf\noindent~\hfill$\bullet\bullet\bullet$\hfill~

\lf\noindent\bf{#1}}

\renewcommand{\bf}[1]{\textbf{#1}}
\renewcommand{\it}[1]{\textit{#1}}

\renewcommand{\sf}[1]{\textsf{#1}}

\renewcommand{\tt}[1]{\texttt{#1}}
\newcommand{\hl}[1]{\bf{\it{#1}}}
\newcommand{\mrm}[1]{\mathrm{#1}}
\newcommand{\mbf}[1]{\mathbf{#1}}
\newcommand{\msf}[1]{\text{\small $\sf{#1}$}}

\newcommand{\cmc}[1]{\mathcal{#1}}
\newcommand{\eus}[1]{\mathscr{#1}}
\newcommand{\euf}[1]{\mathfrak{#1}}
\newcommand{\bb}[1]{\mathbb{#1}}
\newcommand{\msmall}[1]{{\setlength{\arraycolsep}{.6ex}\text{\small$#1$}}}

\newcommand{\mscriptsize}[1]{{\setlength{\arraycolsep}{.3ex}\text{\scriptsize$#1$}}}
\newcommand{\mtiny}[1]{{\setlength{\arraycolsep}{.3ex}\text{\tiny$#1$}}}
\newcommand{\nbd}[1]{$#1$\nobreakdash--}
\newcommand{\ol}[1]{\overline{#1}}
\newcommand{\ul}[1]{\underline{#1}}
\newcommand{\wt}[1]{\widetilde{#1}}

\newcommand{\ve}{\varepsilon}
\newcommand{\vt}{\vartheta}
\newcommand{\vk}{\varkappa}
\newcommand{\vp}{\varphi}

\newcommand{\abs}[1]{\left\lvert#1\right\rvert}
\newcommand{\norm}[1]{\left\lVert#1\right\rVert}

\newcommand{\family}[1]{\left(#1\right)}

\newcommand{\bfam}[1]{\bigl(#1\bigr)}
\newcommand{\Bfam}[1]{\Bigl(#1\Bigr)}
\newcommand{\AB}[1]{\langle#1\rangle}

\newcommand{\BAB}[1]{\Bigl\langle#1\Bigr\rangle}

\newcommand{\BbAB}[1]{\Biggl\langle#1\Biggr\rangle}
\newcommand{\CB}[1]{\{#1\}}
\newcommand{\bCB}[1]{\bigl\{#1\bigr\}}
\newcommand{\BCB}[1]{\Bigl\{#1\Bigr\}}
\newcommand{\SB}[1]{[#1]}

\newcommand{\Matrix}[1]{\begin{pmatrix}#1\end{pmatrix}}
\newcommand{\SMatrix}[1]{\msmall{\Matrix{#1}}}

\newcommand{\sMatrix}[1]{\mscriptsize{\Matrix{#1}}}
\newcommand{\tMatrix}[1]{\mtiny{\Matrix{#1}}}
\newcommand{\rtMatrix}[1]{\raisebox{.3ex}{\tMatrix{#1}}}

\newcommand{\set}[2][]{
    \ifthenelse{\equal{#1}{}}{
        \CB{#2}}{
        \CB{#1~|~#2}}}
\newcommand{\bset}[2][]{
    \ifthenelse{\equal{#1}{}}{
        \bCB{#2}}{
        \bCB{#1~|~#2}}}
\newcommand{\Bset}[2][]{
    \ifthenelse{\equal{#1}{}}{
        \BCB{#2}}{
        \BCB{#1~\big|~#2}}}
\newcommand{\zero}{\CB{0}}

\DeclareMathOperator{\ls}{\normalfont\msf{span}}
\DeclareMathOperator{\cls}{\ol{\ls}}

\DeclareMathOperator{\Tr}{\normalfont\msf{Tr}}
\DeclareMathOperator{\id}{\normalfont\msf{id}}
\DeclareMathOperator{\sid}{\mtiny{\sf{id}}}

\renewcommand{\dim}{\operatorname{\msf{dim}}}

\newcommand{\C}{\bb{C}}

\newcommand{\E}{\bb{E}}

\newcommand{\N}{\bb{N}}

\newcommand{\R}{\bb{R}}

\newcommand{\T}{\bb{T}}

\newcommand{\cA}{\cmc{A}}
\newcommand{\cB}{\cmc{B}}
\newcommand{\cC}{\cmc{C}}

\newcommand{\sB}{\eus{B}}

\newcommand{\sE}{\eus{E}}
\newcommand{\sF}{\eus{F}}

\newcommand{\sK}{\eus{K}}

\newcommand{\sN}{\eus{N}}

\newcommand{\sT}{\eus{T}}

\newcommand{\eK}{\euf{K}}
\newcommand{\eL}{\euf{L}}

\newcommand{\U}{\mbf{1}}
\newcommand{\F}{{\mrm{F}}}

\ifthenelse{\boolean{nosectionequn}}{}{
    \numberwithin{equation}{section}
    }

\ifthenelse{\boolean{article}\or\boolean{letter}\or\boolean{nosectionequn}}{
    \setboolean{nosectionappendix}{true}}{}
\ifthenelse{\boolean{nosectionappendix}}{}{
    \renewcommand{\appendix}{
        \chapter*{\appendixname}
        \addcontentsline{toc}{chapter}{\appendixname}
        \renewcommand{\thesection}{\Alph{section}}
        \setcounter{section}{0}}}

\ifthenelse{\boolean{report}\or\boolean{book}}{
    }{}

\ifthenelse{\boolean{notheorem}}{}{
        \newcommand{\mnname}{Mathematical note.}
        \newcommand{\enname}{End of the note.}
        \newcommand{\definame}{Definition.}
        \newcommand{\propname}{Proposition.}
        \newcommand{\lemname}{Lemma.}
        \newcommand{\exname}{Example.}
        \newcommand{\exername}{Exercise.}
        \newcommand{\remname}{Remark.}
        \newcommand{\obname}{Observation.}
        \newcommand{\thmname}{Theorem.}
        \newcommand{\corname}{Corollary.}
        \newcommand{\proofname}{Proof.}
    \ifthenelse{\boolean{german}}{
        \renewcommand{\mnname}{Mathematische Notiz.}
        \renewcommand{\enname}{Ende der Notiz.}
        \renewcommand{\exname}{Beispiel.}
        \renewcommand{\exername}{Übung.}
        \renewcommand{\remname}{Bemerkung.}
        \renewcommand{\obname}{Beobachtung.}
        \renewcommand{\thmname}{Satz.}
        \renewcommand{\corname}{Korollar.}
        \renewcommand{\proofname}{Beweis.}}{}
    \ifthenelse{\boolean{italian}}{
        \renewcommand{\mnname}{Nota matematica.}
        \renewcommand{\enname}{Fina della nota.}
        \renewcommand{\definame}{Definizione.}
        \renewcommand{\propname}{Proposizione.}
        \renewcommand{\exname}{Esempio.}
        \renewcommand{\exername}{Esercizio.}
        \renewcommand{\remname}{Nota.}
        \renewcommand{\obname}{Osservazione.}
        \renewcommand{\thmname}{Teorema.}
        \renewcommand{\corname}{Corollario.}
        \renewcommand{\proofname}{Dimostrazione.}

       \renewcommand{\appendixname}{Appendice}

       }{}
    \theoremheaderfont{\normalfont\bfseries}
    \theoremstyle{change}
        \theorembodyfont{\rmfamily}
            \newtheorem{emp}{}[section]
                \newcommand{\bemp}[1][]{
                    \begin{emp}\hskip-\labelsep\bf{#1}\hskip\labelsep}
                \newcommand{\eemp}{\end{emp}}
\newtheorem{itemp}[emp]{}
                \newcommand{\bitemp}[1][]{
                    \begin{itemp}\hskip-\labelsep\bf{#1}\hskip\labelsep\normalfont\itshape}
                \newcommand{\eitemp}{\end{itemp}}
            \newtheorem{mn}[emp]{\mnname}
                \newcommand{\bnm}{\begin{mn}~\begin{quotation}\renewcommand{\baselinestretch}{1}\small\noindent\ignorespaces}
                \newcommand{\enm}{\end{quotation}\hfill\bf{\enname}\end{mn}}
            \newtheorem{ex}[emp]{\exname}
                \newcommand{\bex}{\begin{ex}}
                \newcommand{\eex}{\end{ex}}
            \newtheorem{exer}[emp]{\exername}
                \newcommand{\bexer}{\begin{exer}}
                \newcommand{\eexer}{\end{exer}}
            \newtheorem{defi}[emp]{\definame}
                \newcommand{\bdefi}{\begin{defi}}
                \newcommand{\edefi}{\end{defi}}
            \newtheorem{rem}[emp]{\remname}
                \newcommand{\brem}{\begin{rem}}
                \newcommand{\erem}{\end{rem}}
            \newtheorem{ob}[emp]{\obname}
                \newcommand{\bob}{\begin{ob}}
                \newcommand{\eob}{\end{ob}}
        \theorembodyfont{\normalfont\itshape}
            \newtheorem{thm}[emp]{\thmname}
                \newcommand{\bthm}{\begin{thm}}
                \newcommand{\ethm}{\end{thm}}
            \newtheorem{prop}[emp]{\propname}
                \newcommand{\bprop}{\begin{prop}}
                \newcommand{\eprop}{\end{prop}}
            \newtheorem{cor}[emp]{\corname}
                \newcommand{\bcor}{\begin{cor}}
                \newcommand{\ecor}{\end{cor}}
            \newtheorem{lem}[emp]{\lemname}
                \newcommand{\blem}{\begin{lem}}
                \newcommand{\elem}{\end{lem}}
\newenvironment{empn}[1]{\lf\noindent\bf{#1}\ignorespaces\hskip\labelsep}{\lf}
		\newcommand{\bempn}[1]{\begin{empn}{#1}}
		\newcommand{\eempn}{\end{empn}}
		\newcommand{\bitempn}[1]{\begin{empn}{#1}\normalfont\itshape}
		\newcommand{\eitempn}{\end{empn}}
                \newcommand{\bnmn}{\begin{empn}{\mnname}~\begin{quotation}\renewcommand{\baselinestretch}{1}\small\noindent\ignorespaces}
                \newcommand{\enmn}{\end{quotation}\hfill\bf{\enname}\end{empn}}
		\newcommand{\bexn}{\begin{empn}{\exname}}
		\newcommand{\eexn}{\end{empn}}
		\newcommand{\bexern}{\begin{empn}{\exername}}
		\newcommand{\eexern}{\end{empn}}
		\newcommand{\bdefin}{\begin{empn}{\definame}}
		\newcommand{\edefin}{\end{empn}}
		\newcommand{\bremn}{\begin{empn}{\remname}}
		\newcommand{\eremn}{\end{empn}}
		\newcommand{\bobn}{\begin{empn}{\obname}}
		\newcommand{\eobn}{\end{empn}}

\newcommand{\qedsymbol}{~\rule[-0.35mm]{2mm}{2mm}}
    \newcounter{proof}[emp]
    \newenvironment{Proof}[1]{
        \vspace{1ex}
        \renewcommand{\item}[1][\stepcounter{proof}(\roman{proof})]%
            {##1\hskip\labelsep}
        \noindent\textsc{#1\hskip\labelsep}}{
        \nolinebreak\qedsymbol}
    \newcommand{\proof}[1][\proofname]{
        \begin{Proof}{#1}\ignorespaces}
    \newcommand{\qed}{\end{Proof}}
    \newcommand{\noqed}{
        \renewcommand{\qedsymbol}{}
        \end{Proof}}}
    \ifthenelse{\boolean{italian}}{
        \renewcommand{\proofname}{Dimostrazione.}}{}

\usepackage[varg]{txfonts}

\usepackage[matrix,arrow,curve]{xy}

\newcommand{\uodot}{\:\ul{\odot}\:}

\begin{document}

\bibliographystyle{amsalpha}

\title{CP-H-Extendable Maps between Hilbert modules\\and CPH-Semigroups\thanks{AMS 2010 subject classification 46L08, 46L55, 46L53, 60G25}}
\author{}
\author{
~\\
Michael Skeide and K.\ Sumesh
}
\date{October 2012; revised and extended September 2013}

{
\renewcommand{\baselinestretch}{1}
\maketitle


\vfill

\begin{abstract}
\noindent
One may ask which maps between Hilbert modules allow for a completely positive extension to a map acting block-wise between the associated (extended) linking algebras.

In these notes we investigate in particular those \it{CP-extendable maps} where the \nbd{22}cor\-ner of the extension can be chosen to be a homomorphism, the \it{CP-H-extendable maps}. We show that they coincide with the maps considered by Asadi \cite{Asa09}, by Bhat, Ramesh, and Sumesh \cite{BRS12}, and by Skeide \cite{Ske12a}. We also give an intrinsic characterization that generalizes the characterization by Abbaspour and Skeide \cite{AbSk07} of homomorphicly extendable maps as those which are  ternary homomorphisms.

For general \it{strictly} CP-extendable maps we give a factorization theorem that generalizes those of Asadi, of Bhat, Ramesh, and Sumesh, and of Skeide for CP-H-extendable maps. This theorem may be viewed as a unification of the representation theory of the algebra of adjointable operators and the \it{KSGNS-construction}.

As an application, we examine semigroups of CP-H-extendable maps, so-called \it{CPH-semigroups}, and illustrate their relation with a sort of generalized dilation of CP-semi\-groups, \it{CPH-dilations}.
\end{abstract}

}



	\newpage

\section{Introduction}

Let $\tau\colon\cB\rightarrow\cC$ be a linear map between \nbd{C^*}algebras $\cB$ and $\cC$. A \hl{\nbd{\tau}map} is a map $T\colon E\rightarrow F$ from a Hilbert \nbd{\cB}module $E$ to a Hilbert \nbd{\cC}module $F$ such that
\beqn{\label{*}\tag{$*$}
\AB{T(x),T(x')}
~=~
\tau(\AB{x,x'}).
}\eeqn
After several publications about \nbd{\tau}maps where $\tau$ was required to be a homomorphism (for instance, Bakic and Guljas \cite{BaGu02}, Skeide \cite{Ske06c}, Abbaspour Tabadkan and Skeide \cite{AbSk07}), and others where $\tau$ was required to be just a CP-map (for instance, Asadi \cite{Asa09}, Bhat, Ramesh, and Sumesh \cite{BRS12}, Skeide \cite{Ske12a}), we think it is now time to determine the general structure of \nbd{\tau}maps. We also think it is time to, finally, give some idea what \nbd{\tau}maps might be good for. While we succeed completely with our first task for bounded $\tau$, we hope that our small application in Section \ref{CPHdilSEC} that establishes a connection with dilations of CP-semigroups and product systems can, at least, view perspectives for concrete applications in the future.

\bulletline
If $T$ fulfills \eqref{*} for some linear map $\tau$, then $T$ is linear. (Examine $\abs{T(x+\lambda x')-T(x)-\lambda T(x')}^2$.) Furthermore, if $\tau$ is bounded, then, obviously, $T$ is bounded with norm $\norm{T}\le\sqrt{\norm{\tau}}$. As easily, one checks that the \hl{inflation} $T_n\colon M_n(E)\rightarrow M_n(F)$ of $T$ (that is, $T$ acting element-wise on the matrix) is a \nbd{\tau_n}map for the inflation $\tau_n\colon M_n(\cB)\rightarrow M_n(\cC)$ of $\tau$. (Recall that $M_n(E)$ is a Hilbert \nbd{M_n(\cB)}module with inner product $\bfam{\AB{X,Y}}_{i,j}:=\sum_k\AB{x_{ki},y_{kj}}$.) Therefore, $\norm{T_n}\le\sqrt{\norm{\tau_n}}$.

A map $\tau$ fulfilling \eqref{*} (and, therefore, also $\tau_n$) ``looks'' positive. (In fact, at least positive elements of the form $\AB{x,x}$ are sent to the positive elements $\AB{T(x),T(x)}$.) More precisely, it looks positive on the ideal $\ls\AB{E,E}$. It is not difficult to show (see Lemma \ref{tauCPlem}) that bounded $\tau$ is, actually, positive on the \hl{range ideal} $\cB_E:=\cls\AB{E,E}$ of $E$. Since the same is true also for $\tau_n$, we see that $\tau$ is \hl{completely positive} (or \hl{CP}) on $\cB_E$. Recall that for CP-maps $\tau$ we have $\norm{\tau_n}=\norm{\tau}$.

We arrive at our first new result.

\bthm\label{gentauthm}
Let $T\colon E\rightarrow F$ be a map from a \hl{full} Hilbert \nbd{\cB}module $E$ (that is, $\cB_E=\cB$) to a Hilbert \nbd{\cC}module $F$, and let $\tau\colon\cB\rightarrow\cC$ be a bounded linear map. If $T$ a \nbd{\tau}map, then $\tau$ is completely positive. Moreover, $T$ is linear and completely bounded with \hl{CB-norm} $\norm{T}_{cb}:=\sup_n\norm{T_n}=\sqrt{\norm{\tau}}$.
\ethm

The second missing part (apart from Lemma \ref{tauCPlem}), namely, that the CB-norm $\norm{T}_{cb}$ actually reaches its bound $\sqrt{\norm{\tau}}$, we prove in Lemma \ref{maxnlem}.

It is, in general, not true that $\norm{T}_{cb}=\norm{T}$, not even if $\cB$ and $\cC$ are unital.%
\footnote{
This contradicts the proposition in Asadi \cite{Asa09}.}
It is true, if $E$ has a \hl{unit vector} $\xi$ (that is, $\AB{\xi,\xi}=\U$); see Observation \ref{uniob}.

\bex\label{counterex}
Let $H\ne\zero$ be a Hilbert space with ONB $\bfam{e_i}_{i\in I}$. For $E$ we choose the full Hilbert \nbd{\sK(H)}module $H^*$ (with inner product $\AB{x'^*,x^*}:=x'x^*$). For $F$ we choose $H$. So, $\cB=\sK(H)$ and $\cC=\C$. Let $T$ be the \hl{transpose map} with respect to the ONB. That is, $T$ sends the ``row vector'' $x^t=\sum_ix_ie_i^*$ in $E$ to the ``column vector'' $x=(x^t)^t=\sum_i x_ie_i$ in $F$. Of course, $\norm{T}=1$.

A linear map $\tau\colon\sK(H)\rightarrow\C$ turning $T$ into a \nbd{\tau}map, would send $e_ie_j^*$ to $\tau(\AB{e_i^*,e_j^*})=\AB{T(e_i^*),T(e_j^*)}=\AB{e_i,e_j}=\delta_{i,j}$. So, on the finite-rank operators $\sF(H)$ the map $\tau$ is bound to be the (non-normalized) \hl{trace} $\Tr:=\sum_i\AB{e_i,\bullet e_i}$. Recall that $\norm{\Tr}=\dim H$. This shows several things:
\begin{enumerate}
\item\label{E1}
Suppose $H$ is infinite-dimensional. Then $\tau$ cannot be bounded. Since positive maps are bounded, there cannot be whatsoever positive map $\tau$ turning $T$ into a \nbd{\tau}map. (Of course, we can extend $\tau=\Tr$ by brute-force linear algebra from $\sF(E)$ to $\sK(E)$, so that $T$ is still a \nbd{\tau}map with unbounded and non-positive $\tau$.)

\item\label{E2}
Suppose $H$ is \nbd{n}dimensional (so that, in particular, $\sK(H)=M_n$ is unital). The column vector ${X^*}^n$ in ${H^*}^n$ with entries $e_1^*,\ldots,e_n^*$ has square modulus $\AB{{X^*}^n,{X^*}^n}=\sum_{i=1}^ne_ie_i^*$. So, $\norm{{X^*}^n}=\sqrt{\norm{\sum_{i=1}^ne_ie_i^*}}=1$. However, the norm of the column vector $Y^n$ with entries $T(e_1^*)=e_1,\ldots,T(e_n^*)=e_n$ is $\sqrt{\sum_{i=1}^n\AB{e_i,e_i}}=\sqrt{n}$. Since $M_n(H^*)\supset M_{n,1}(H^*)={H^*}^n$, we find $\norm{T}_{cb}\ge\norm{T_n}\ge\sqrt{n}$. On the other hand, by the discussion preceding Theorem \ref{gentauthm}, $\norm{T}_{cb}\le\sqrt{\norm{\tau}}=\sqrt{n}$. Therefore, $\norm{T}_{cb}=\sqrt{n}$, while $\norm{T}=1\ne\norm{T}_{cb}$ for $n\ge2$.
\end{enumerate}
\eex

\noindent
Whenever $\cB_E$ is unital, $\tau$ is bounded (and, therefore, completely bounded) on $\cB_E\overset{(!)}{=}\ls\AB{E,E}$; see again Observation \ref{uniob}.

To summarize: If $E$ is full and if $\tau$ is bounded, then CP is automatic. And if $E$ is full over a unital \nbd{C^*}algebra, then we have not even to require that $\tau$ is bounded. On the other hand, some of the questions we wish to tackle, have nice answers for CP-maps $\tau$, even if $E$ is not full. And \nbd{\tau}maps $T$ (into the Hilbert \nbd{\sB(G)}module $F=\sB(G,H)$) for completely positive $\tau$ (into $\cC=\sB(G)$) is also what Asadi started analyzing in \cite{Asa09}. So, after these considerations, for the rest of these notes $\tau$ will \bf{always} be a CP-map. 

\bulletline
A basic task of these notes is to characterize \nbd{\tau}maps for CP-maps $\tau$. More precisely, we wish to find criteria that tell us when a map $T\colon E\rightarrow F$ is a \nbd{\tau}map for some CP-map $\tau$ \bf{without} knowing $\tau$, just by looking at $T$.

The case when a possible $\tau$ is required to be a homomorphism has been resolved by Abbaspour Tabadkan and Skeide \cite{AbSk07}. (In this case, $T$ has been called \hl{\nbd{\tau}homo\-morphism} in Bakic and Guljas \cite{BaGu02} or \hl{\nbd{\tau}isometry}.) For full $E$, \cite[Theorem 2.1]{AbSk07} asserts: $T$ is a \nbd{\tau}isometry for some homomorphism $\tau$ if and only if $T$ is linear and fulfills
\beqn{
T(x\AB{y,z})
~=~
T(x)\AB{T(y),T(z)},
}\eeqn
that is, if $T$ is a \hl{ternary homomorphism}.%
\footnote{\label{terFN}
We should emphasize that, unlike stated in \cite{AbSk07}, linearity of $T$ \bf{cannot} be dropped. The map $T\colon E\rightarrow\C$ defined as $T(x)=1$ is a counter example. Indeed, without linearity, the map $\tau=\vp$ defined in the proof of \cite[Theorem 2.1]{AbSk07} is a well-defined multiplicative \nbd{*}map; but it may fail to be linear.
}
(Ternary homomorphisms into $\sB(G,H)$ ($G$ and $H$ Hilbert spaces) occurred under the name \hl{representation} of a Hilbert module (and the unnecessary hypothesis of complete boundedness) in Skeide \cite{Ske00b}.) Another equivalent criterion is that $T$ extends as a homomorphism acting block-wise between the linking algebras of $E$ and of $F$. (This follows simply by applying \cite[Theorem 2.1]{AbSk07} also to the ternary homomorphism $T^*\colon E^*\rightarrow F^*$ from the full Hilbert \nbd{\sK(E)}module $E^*$ (with inner product $\AB{x'^*,x^*}:=x'x^*$) to the full Hilbert \nbd{\sK(F)}module $F^*$ defined as $T^*(x^*):=T(x)^*$, resulting in a homomorphism $\vt\colon\sK(E)\rightarrow\sK(F)$ so that the block-wise map
\beqn{
\SMatrix{\tau&T^*\\T&\vt}
\colon
\SMatrix{\cB&E^*\\E&\sK(E)}
~\longrightarrow~
\SMatrix{\cC&F^*\\F&\sK(F)}
}\eeqn
is a homomorphism.) We would call such maps \hl{H-extendable}.

It is always a good idea to look at properties of Hilbert modules in terms of properties of their linking algebras. (For instance, Skeide \cite{Ske00b} defined a Hilbert module $E$ over a von Neumann algebra to be a \it{von Neumann module} if its extended linking algebra is a von Neumann algebra in a canonically associated representation. This happens if and only if $E$ is self-dual, that is, if $E$ is a \nbd{W^*}module.) Likewise, it is a good idea to look at properties of maps between Hilbert modules in terms of how they may be extended to block-wise maps between their linking algebras. (For instance, many maps between von Neumann modules are \nbd{\sigma}weakly continuous if and only if they allow for a normal (that is, order continuous) block-wise extension to a map between the linking algebras.) In addition to the usual \hl{linking algebra} $\rtMatrix{\cB&E^*\\E&\sK(E)}=\sK\rtMatrix{\cB\\E}$ of a Hilbert \nbd{\cB}module $E$, it is sometimes useful to look at the \hl{reduced linking algebra} $\rtMatrix{\cB_E&E^*\\E&\sK(E)}$ or at the \hl{extended linking algebra} $\rtMatrix{\cB&E^*\\E&\sB^a(E)}$. It would be tempting to see if \nbd{\tau}maps are precisely the \hl{CP-extendable maps}, that is, maps that allow for some block-wise CP-extension between some sort of linking algebras. Unfortunately, this is not so: There are more CP-extendable maps than \nbd{\tau}maps; see Section \ref{CPrSEC}. We, therefore, strongly object to use the name \it{CP-maps} between Hilbert modules as meaning \nbd{\tau}maps, which was proposed recently by several authors; see, for instance, Heo and Ji \cite{HeJi11}, or Joita \cite{Joi12}.

But if CP-extendable is not the right condition, what is the right condition? And what is the right ``intrinsic condition'' replacing the ternary condition for \nbd{\tau}isometries? As a main result of these notes, in Section \ref{proofSEC} we prove the following theorem.

\newpage

\bthm\label{mainthm}
Let $E$ be a full Hilbert \nbd{\cB}module and let $F$ be a Hilbert \nbd{\cC}module. Let $T\colon E\rightarrow F$ be a linear map and denote $F_T:=\cls T(E)\cC$. Then the following conditions are equivalent:
\begin{enumerate}
\item\label{1}
There exists a (unique) CP-map $\tau\colon\cB\rightarrow\cC$ such that $T$ is a \nbd{\tau}map.

\item\label{2}
$T$ extends to a block-wise CP-map $\sT=\rtMatrix{\tau&T^*\\T&\vt}\colon\rtMatrix{\cB&E^*\\E&\sB^a(E)}\rightarrow\rtMatrix{\cC&F_T^*\\F_T&\sB^a(F_T)}$ where $\vt$ is a homomorphism, that is, $T$ is a \hl{CP-H-extendable map} .

\item\label{3}
$T$ is a completely bounded map and $F_T$ can be turned into a \nbd{\sB^a(E)}\nbd{\cC} correspondence in such a way that $T$ is left \nbd{\sB^a(E)}linear.

\item\label{4}
$T$ is a completely bounded map fulfilling
\beqn{\label{**}\tag{$**$}
\BAB{T(y),T\bfam{x\AB{x',y'}}}
~=~
\BAB{T\bfam{x'\AB{x,y}},T(y')}.
}\eeqn
\end{enumerate}
\ethm
A more readable version of \eqref{**} is
\beqn{
\AB{T(y),T(xx'^*y')}
~=~
\AB{T(x'x^*y),T(y')}.
}\eeqn
This \hl{quaternary condition} is the intrinsic condition we were seeking, and which generalizes the ternary condition guaranteeing that $T$ is \nbd{\tau}isometry.

\bob\label{discob}
\begin{enumerate}
\item\label{O1}
The homomorphism $\vt$ in \eqref{2} coincides with the left action in \eqref{3}; see the proof of \eqref{2} $\Rightarrow$ \eqref{3} in Section \ref{proofSEC}.

\item\label{O2}
Since the set $T(E)$ generates the Hilbert \nbd{\cC}module $F_T$, the left action in \eqref{3} (and, consequently, also $\vt$ in \eqref{2}) is uniquely determined by $(xy^*)T(z)=T(xy^*z)$. In fact, this formula shows that the finite-rank operators $\sF(E)$ act nondegenerately on $F_T$, so there is a unique extension to all of $\sB^a(E)$. Moreover, this unique extension is strict and unital; see the proof of Lemma \ref{strGNSlem}.

\item\label{O3}
It is routine to show that \eqref{**} well-defines a nondegenerate action of $\sF(E)$. So, the same argument also shows that \eqref{3} and \eqref{4} are equivalent.

\item\label{O4}
Clearly, Example \ref{counterex}\eqref{E1} shows that the condition on $T$ to be completely bounded in \eqref{3} and \eqref{4}, may not be dropped. However, if $E$ is full over a unital \nbd{C^*}algebra, then $T$ just linear is sufficient; see again Observation \ref{uniob} .
\end{enumerate}
\eob

\brem
It should be noted that the CP-map $\tau$ in \eqref{2} need not coincide with the map $\tau$ in \eqref{1} making $T$ a \nbd{\tau}map. (Just add an arbitrary CP-map $\cB\rightarrow\cC$ to the latter.) Likewise, having a CP-extension $\sT$ with a non-homomorphic \nbd{22}corner $\vt$ does not necessarily mean that it is not possible to get a CP-H-extension by modifying $\vt$.
\erem

\brem
Unlike for \nbd{\tau}isometries, for more general \nbd{\tau}maps the homomorphism $\vt$ in \eqref{2} will only rarely map the compacts $\sK(E)$ into the compacts $\sK(F_T)$. So, in \eqref{2} it is forced that we pass to the extended linking algebras. Also considerations about the strict topology cannot be avoided completely.
\erem

\brem\label{unirem}
We already know that a \nbd{\tau}map $T$ is linear, so linearity of $T$ may be dropped from \eqref{1}. We know from the example in Footnote \ref{terFN} that linearity may not be dropped from \eqref{4}, not even if $T$ fulfills the stronger ternary condition. Linearity may be dropped from \eqref{3}, if $E$ contains a unit vector $\xi$ (or, more generally, a direct summand $\cB$), for in that case we have $T(x)=T(x\xi^*\xi)=(x\xi^*)T(\xi)$, which is linear in $x$. However, unlike in Observation \ref{discob}\eqref{O4}, we were not able to save the statement for unital $\cB$ without a unit vector.
\erem

The property in \eqref{3} is almost visible from a glance at \eqref{*}. In fact, we try to assign a value $\AB{T(x),T(x')}$ to the element $\AB{x,x'}\in\cB_E=E^*\odot E$. (Here $E^*$ is the \hl{dual} Hilbert \nbd{\sB^a(E)}module of $E$ with inner product $\AB{x'^*,x^*}:=x'x^*$, and the tensor product is over the canonical left action of $\sB^a(E)$ on $E$.) It is clear that the map $(x,x')\mapsto\AB{T(x),T(x')}$ has to be balanced over $\sB^a(E)$ if there should exist $\tau$ fulfilling \eqref{*}. And if there was a suitable left action of $\sB^a(E)$ on $F_T$, then we would be concerned with the map $\tau:=T^*\odot T$. People knowing the module Haagerup tensor product of operator modules and Blecher's result \cite[Theorem 4.3]{Ble97} that the Haagerup tensor product is (completely) isometrically isomorphic to the tensor product of correspondences, can already smell that everything is fine. We shall give a direct proof in Section \ref{proofSEC}. Actually, our method will provide us with a quick proof of Blecher's result.

\bulletline
We have seen in Theorem \ref{mainthm} that the Hilbert submodule $F_T$ of $F$ generated by $T(E)$ plays a distinguished role. (If $T$ is a \nbd{\tau}isometry, then $T(E)$ is already a closed \nbd{\tau(\cB)}submodule of $F$.) It is natural to ask to what extent the condition in \eqref{2} can be satisfied if we write $F$ instead of $F_T$. In developing semigroup versions in Sections \ref{CPHsgSEC} and \ref{CPHdilSEC}, this situation becomes so important that we prefer to use the the acronym \it{CPH} for that case, and leave for the equivalent of \nbd{\tau}maps the rather contorted term \it{CP-H-extendable}:

\bdefi
A \hl{CPH-map} from $E$ to $F$ is a map that extends as a block-wise CP-map between the extended linking algebras of $E$ and of $F$ such that the \nbd{22}corner is a homomorphism. A CPH-map is \hl{strictly CPH} if that homomorphism can be chosen strict. A (strictly) CPH-map is a (\hl{strictly}) \hl{CPH$_0$-map} if the homomorphism $\vt$ can be chosen unital.
\edefi

CPH-maps are CP-H-extendable (Lemma \ref{CP(-)Hlem}). If $F_T$ is complemented in $F$, then $T$ is a CPH-map if and only if it is CP-H-extendable. (In that case, $\sB^a(F_T)$ is a corner of $\sB^a(F)$, so that $\vt$ may be considered a map into $\sB^a(F)$.) But this condition is not at all necessary, nor natural; see Observation \ref{notETob}.

Despite the fact that there are fewer CPH-maps than CP-H-extendable maps, looking at CPH-maps is particularly crucial if we wish to look at semigroups of CP-H-extendable maps $T_t$ on $E$. Obviously, for full $E$, the associated CP-maps $\tau_t$ form a CP-semigroup. But the same question for the homomorphisms $\vt_t$, \it{a priori}, has no meaning.  The extensions $\vt_t$ map $\sB^a(E)$ into $\sB^a(E_{T_t})$, not into $\sB^a(E)$. And if $E_{T_t}$ is not complemented in $E$, then it is not possible to interpret $\sB^a(E_{T_t})$ as a subset of $\sB^a(E)$, to which $\vt_s$ could be applied in order to make sense out of $\vt_s\circ\vt_t$.

In Section \ref{CPHsgSEC} we study such \hl{CPH-semigroups}, and examine how the results of the first sections may be generalized or reformulated. These results depend essentially on the theory of tensor product systems of correspondences initiated Bhat and Skeide \cite{BhSk00} (following Arveson \cite{Arv89} for Hilbert spaces), which, in our case, have to replace the GNS-construction for a single CP-map $\tau$. See the introduction to Section \ref{CPHsgSEC}, in particular after Observation \ref{unistriob}, for more details.

In the speculative Section \ref{CPHdilSEC} we introduce the new concept of \hl{CPH-dilation} of a CP-map or a CP-semigroup. It generalizes the concept of \it{weak dilation} and is intimately related to CPH-maps or CPH-semigroups. In the end, we comment on some relations with (completely positive definite) \hl{CPD-kernels} and with \hl{Morita equivalence}. If CPH-dilations  can be considered an interesting concept, and if, as demonstrated, understanding CPH-dilations is the same understanding CPH-maps and CPH-semigroups, then Section \ref{CPHdilSEC} shows the road to what might be the first application of CPH-maps.

\bulletline
We wish to underline that \bf{all} results in these can be formulated for von Neumann algebras, von Neumann modules (or \nbd{W^*}modules), and von Neumann correspondences (or \nbd{W^*}cor\-re\-spond\-ences), replacing  also the tensor product of \nbd{C^*}correspondences with that of von Neumann correspondences, replacing \it{full} with \it{strongly full} ,and adding to all maps between von Neumann objects the word \it{normal} (or \nbd{\sigma}weak). We do not give any detail, because the proofs either generalize word by word or are simple adaptations of the \nbd{C^*}proofs. We emphasize, however, that all problems regarding adjointability of maps or complementability of $F_T$ in $F$ disappear. Therefore, for a map between von Neumann modules. Likewise, every normal CP-map from the von Neumann algebra $\ol{\cB}_E^s$ extends to a CP-from $\cB$. CPH and CP-H-extendable for von Neumann modules (or \nbd{W^*}modules) is the same thing and they do no longer depend on (strong) fullness.

\newpage

\section{Proof of Theorem \ref{mainthm}}\label{proofSEC}

Equivalence of \ref{3} and \ref{4} has already been dealt with in Observation \ref{discob}\eqref{O2} and \eqref{O3}. For the remaining steps we shall follow the order \ref{1} $\Rightarrow$ \ref{2} $\Rightarrow$ \ref{3} $\Rightarrow$ \ref{1}. Since we also wish to make comments on the mechanisms of some steps or how parts of the proof are applicable in more general situations, we put each of the steps into an own subsection and indicate by ``\!\qedsymbol'' where the part specific to Theorem \ref{mainthm} ends.

In Section \ref{CPrSEC}, we present an alternative direct proof of \ref{2} $\Rightarrow$ \ref{1}, which avoids using arguments originating in operator spaces as involved in the proof \ref{3} $\Rightarrow$ \ref{1}.

\subsection*{Proof \ref{1} $\Longrightarrow$ \ref{2}}

We first consider the case where $\cB$ and $\cC$ are unital, but without requiring that $E$ is full. So let $\tau\colon\cB\rightarrow\cC$ be a CP-map between unital \nbd{C^*}algebras, and let $T\colon E\rightarrow F$ be a \nbd{\tau}map from an arbitrary Hilbert \nbd{\cB}module $E$ to a Hilbert \nbd{\cC}module $F$.

Since $\cB$ and $\cC$ are unital, by Paschke's \hl{GNS-construction} \cite{Pas73} for $\tau$, we get a pair $(\sF,\zeta)$ consisting of \hl{GNS-correspondence} $\sF$ from $\cB$ to $\cC$ and \hl{cyclic vector} $\zeta$ in $\sF$ such that
\baln{
\AB{\zeta,\bullet\zeta}
&
~=~
\tau,
&
\cls\cB\zeta\cC
&
~=~
\sF.
}\ealn
One easily verifies that the map
\beqn{
x\odot\zeta
~\longmapsto~
T(x)
}\eeqn
defines an isometry $v\colon E\odot\sF\rightarrow F$. (It maps $x\odot(b\zeta c)=((xb)\odot\zeta)c$ to $T(xb)c$.) In other words, $T$ factors as $T=v(\bullet\odot\zeta)$. (We just have reproduced the simple proof of the ``only if'' direction of the theorem in Skeide \cite{Ske12a}.)

Now, $v$ is obviously a unitary onto $F_T:=\cls T(E)\cC$. So $\vt:=v(\bullet\id_\sF)v^*$ defines a (unital and strict) homomorphism $\sB^a(E)\rightarrow\sB^a(F_T)$. Identifying $\sF$ with $\sB^a(\cC,\sF)$ via $y\colon c\mapsto yc$ and identifying $\cB\odot\sF$ with $\sF$ via $b\odot y\mapsto by$, we may define a map
\beqn{
\Xi
~:=~
\sMatrix{\zeta\,&\\&\,v^*}
~\in~
\sB^a\sMatrix{\sMatrix{\cC\\F_T}\,,\,\sMatrix{\cB\\E}\odot\sF}.
}\eeqn
Obviously, the map $\sT:=\Xi^*(\bullet\odot\id_\sF)\Xi$ from the extended linking algebra of $E$ into the extended linking algebra of $F_T$ is completely positive. One easily verifies that
\beqn{
\sT
~=~
\sMatrix{\tau&T^*\\T&\vt},
}\eeqn
where $T^*(x^*):=T(x)^*$. This proves \ref{1} $\Rightarrow$ \ref{2} for unital \nbd{C^*}algebras but not necessarily full $E$.

Now suppose $\cB$ is not necessarily unital. (Nonunital $\cC$ may always be ``repaired'' by appropriate use of approximate units.) The following is folklore.

\blem\label{unizelem}
If $\tau\colon\cB\rightarrow\cC$ is a CP-map, then the map $\wt{\tau}\colon\wt{\cB}\rightarrow\wt{\cC}$ between the unitalizations of $\cB$ and $\cC$, defined by
\baln{
\wt{\tau}\upharpoonright\cB
&
~:=~
\tau,
&
\wt{\tau}(\wt{\U})
&
~:=~
\norm{\tau}\wt{\U},
}\ealn
is a CP-map, too.
\elem

\proof
Denote by $\delta\colon\wt{\cB}\rightarrow\C$ the unique character vanishing on $\cB$, and choose a contractive approximate unit $\bfam{u_\lambda}_{\lambda\in\Lambda}$ for $\cB$. Then the maps
\beqn{
\tau_\lambda
~:=~
\tau(u_\lambda^*\bullet u_\lambda)+\Bfam{\norm{\tau}\wt{\U}-\tau(u_\lambda^*u_\lambda)}\delta
}\eeqn
are CP-maps (as sum of CP-maps) and converge pointwise to $\wt{\tau}$. Therefore, $\wt{\tau}$ is a CP-map, too.\qed

\lf
Now, $E$ and $F$ are also modules over the unitalizations, and $T$ is a \nbd{\wt{\tau}}map, too. Since in the first part $E$ was not required full, we may apply the result and get a CP-map $\wt{\sT}$ that, obviously, restricts to the desired CP-map $\sT$. This concludes the proof \ref{1} $\Rightarrow$ \ref{2}.\qedsymbol

\bob
Obviously, the proof shows that the conclusion \ref{1} $\Rightarrow$ \ref{2} holds in general, even if $E$ is not full: All \nbd{\tau}maps are CP-H-extendable.
\eob

\bob \label{thmrevob}
Adding the obvious statement that for each \nbd{\cB}\nbd{\cC}correspondence $\sF$ and for each vector $\zeta\in\sF$, an isometry $v\colon E\odot\sF\rightarrow F$ gives rise to a \nbd{\tau}map $T:=v(\bullet\odot\zeta)$ for the CP-map $\tau:=\AB{\zeta,\bullet\zeta}$, we also get the ``if'' direction of the theorem in \cite{Ske12a}. For this it is not necessary that $\sF$ is the GNS-correspondence of $\tau$. This observation provides us with many CPH-maps. It also plays a role Section \ref{CPHsgSEC}.
\eob

\brem
The theorem in Skeide \cite{Ske12a} is the last and most general version of a result, first, stated by Asadi \cite{Asa09} for unital CP-maps into $\cC=\sB(G)$ and $T$ mapping into $F=\sB(G,H)$ ($G$ and $H$ Hilbert spaces) under the extra condition that $T(\xi)T(\xi)^*=\id_F$ for some $\xi\in E$ and, then, proved by Bhat, Ramesh, and Sumesh \cite{BRS12} (without the extra condition and for $\cB$ still unital, but $\tau$ not necessarily unital).
\erem

\subsection*{Proof \ref{2} $\Longrightarrow$ \ref{3}}

Let $T\colon E\rightarrow F$ be a map from a Hilbert \nbd{\cB}module $E$ to a Hilbert \nbd{\cC}module $F$. Define the map $T^*\colon x^*\mapsto T(x)^*$, and put $F_T:=\cls T(E)\cC$. Suppose we find a CP-map  $\tau\colon \cB\rightarrow\cC$ and a homomorphism $\vt\colon \sB^a(E)\rightarrow \sB^a(F_T)$ such that $\sT:=\rtMatrix{\tau&T^*\\T&\vt}\colon\rtMatrix{\cB&E^*\\E&\sB^a(E)}\rightarrow\rtMatrix{\cC&F_T^*\\F_T&\sB^a(F_T)}$ is a CP-map. Then, in particular, $T$ is a CB-map.

\blem\label{bilcp}
Let $S\colon \cB\rightarrow\cC$ be a CP-map between \nbd{C^*}algebras $\cB$ and $\cC$. Suppose $\cA\subset\cB$ is a \nbd{C^*}subalgebra of $\cB$ with unit $\U_{\cA}$ such that the restriction $\vt:=S\upharpoonright\cA$ of $S$ to $\cA$ is a homomorphism. Then
\baln{
S(ba)
&
~=~
S(b\U_\cA)\vt(a),
&
S(ab)
&
~=~
\vt(a)S(\U_\cA b)
&
}\ealn
for all $b\in\cB$ and $a\in\cA$.
\elem

\proof
Assume that $\cB$ and $\cC$ are unital. (Otherwise, unitalize as explained in Lemma \ref{unizelem} and observe that also the unitalization $\wt{S}$ fulfills the hypotheses for $\cA\subset\wt{\cB}$ with the same $\vt$. If the statement is true for $\wt{S}$, then so it is for $S=\wt{S}\upharpoonright\cB$.)

Let $(\sF,\zeta)$ denote the GNS construction for $S$. By the stated properties, one easily verifies that $\abs{a\zeta-\U_{\cA}\zeta\vt(a)}^2=0$, so, $a\zeta=\U_{\cA}\zeta\vt(a)$ for all $a\in\cA$. The first equation of the lemma follows by computing $S(ba)=\AB{\zeta,ba\zeta}$, and the second by taking its adjoint.\qed

\lf
By applying Lemma \ref{bilcp} to the CP-map $\sT\colon\rtMatrix{\cB&E^*\\E&\sB^a(E)}\rightarrow\rtMatrix{\cC&F_T^*\\F_T&\sB^a(F_T)}$ with the subalgebra $\cA=\rtMatrix{0&0\\0&\sB^a(E)}\ni\rtMatrix{0&0\\0&\sid_E}=\U_\cA$, we get
\beqn{
\sMatrix{0&0\\T(ax)&0}
~=~
\sT\sMatrix{\sMatrix{0&0\\0&a}\sMatrix{0&0\\x&0}}
~=~
\sMatrix{0&0\\0&\vt(a)}\sT\sMatrix{0&0\\x&0}
~=~
\sMatrix{0&0\\\vt(a)T(x)&0},
}\eeqn
thus $T(ax)=\vt(a)T(x)$ for all $x\in E$ and $a\in\sB^a(E)$. This proves \ref{2} $\Rightarrow$ \ref{3}.\qedsymbol

\bob
Also here we did not require that $E$ is full. So \ref{2} $\Rightarrow$ \ref{3} is true for all CP-H-extendable maps.
\eob

Effectively, for the conclusion $T(ax)=\vt(a)T(a)$, we did not even need that $\sT$ maps into the linking algebra of $F_T$. The conclusion remains true for all CPH-maps, so that for a CPH-map the subspace $F_T$ of $F$ reduces $\vt$.

\bcor\label{CP(-)Hlem}
A CPH-map $T\colon E\rightarrow F$ is CP-H-extendable.
\ecor

For full $E$, this also follows via CPH $\Rightarrow$ \ref{3} $\Rightarrow$ \ref{1} $\Rightarrow$ \ref{2}, as soon as we have completed the step \ref{3} $\Rightarrow$ \ref{1}.

\subsection*{Proof \ref{3} $\Longrightarrow$ \ref{1}}

Given $T$ and a left action of $\sB^a(E)$ on $F_T$ such that $aT(x)=T(ax)$, our scope is to define $\tau$ by \eqref{*}. So, in this part it is essential that $E$ is full. Our job will be to show that the hypotheses of \ref{3}, which showed already to be necessary, are also sufficient.

As mentioned in the introduction, in the case $\cB=\cB_E=E^*\odot E$, the map $\tau$, if it exists, appears to be the map
\beqn{
\cB
~=~
E^*\odot E
~\xrightarrow{~T^*\odot T~}~
F^*\odot F
~=~
\cC_F
~\subset~
\cC.
}\eeqn
Note that, actually, $T^*\odot T$ maps into $F_T^*\odot F_T\subset F^*\odot F$. And if $F$ is a correspondence making $T$ left \nbd{\sB^a(E)}linear, then, by definition of left \nbd{\sB^a(E)}linear, $F_T$ is a correspondence making $T$ left \nbd{\sB^a(E)}linear, too. (Also strictness does not play any role here.) So, it does note really matter if we require the property in \ref{3} for $F_T$ or for $F$, because the latter implies the former. So, let $F$ be a \nbd{\sB^a(E)}\nbd{\cC}correspondence such that $T$ is left \nbd{\sB^a(E)}linear. Likewise, $T^*:=*\circ T\circ*$ is a right \nbd{\sB^a(E)}linear map for the corresponding \nbd{\sB^a(E)}module structures of $E^*$ and $F^*$. So, $T^*\odot T$, indeed, defines a linear map from the algebraic tensor product $E^*\uodot E$ over $\sB^a(E)$ into $F^*\odot F$. And by Lance \cite[Proposition 4.5]{Lan95}, we have $E^*\uodot E=\ls\AB{E,E}$ as subset of $E^*\odot E=\cB$.

Once $\tau\colon E^*\uodot E\rightarrow\cC$ is bounded (for the norm of the internal tensor product $E^*\odot E$ on $E^*\uodot E\subset E^*\odot E$), Theorem \ref{gentauthm} asserts that the extension to $\cB=E^*\odot E$ is completely positive. Recall that we still have to add the following missing piece to the proof of that theorem:

\blem\label{tauCPlem}
Let $\tau\colon\cB\rightarrow\cC$ be a bounded linear map fulfilling \eqref{*} for some map $T\colon E\rightarrow F$. Then $\tau$ is positive on $\cB_E$.
\elem

\proof
We already said that $T$ being a \nbd{\tau}map, also $T_n$ is a \nbd{\tau_n}map. Similarly, $T^n\colon E^n\rightarrow F^n$ is a \nbd{\tau}map itself. Let us choose a bounded approximate unit $\bfam{u_\lambda}_{\lambda\in\Lambda}$ for $\cB_E$ consisting of elements $u_\lambda=\sum_{i=1}^{n_\lambda}\AB{x_i^\lambda,y_i^\lambda}\in\cB_E$. Defining the elements $X_\lambda\in E^{n_\lambda}$ with entries $x_i^\lambda$ and, similarly, $Y_\lambda$, we get $u_\lambda=\AB{X_\lambda,Y_\lambda}$. For any positive element $bb^*$ in $\cB_E$, denote by $a_\lambda\in\sK(E^{n_\lambda})$ the positive square root of the rank-one operator $X_\lambda bb^*X_\lambda^*=(X_\lambda b)(X_\lambda b)^*$. Then
\beqn{
\tau(u_\lambda^*bb^*u_\lambda)
~=~
\tau(\AB{a_\lambda Y_\lambda,a_\lambda Y_\lambda})
~=~
\AB{T^{n_\lambda}(a_\lambda Y_\lambda),T^{n_\lambda}(a_\lambda Y_\lambda)}
~\ge~
0.
}\eeqn
Since $u_\lambda^*bb^*u_\lambda\to bb^*$ in norm, and since $\tau$ is bounded, we get $\tau(bb^*)\ge0$.
\qed

\lf
So it remains to show that $\tau$ is bounded on $E^*\uodot E$. Care is in place, however, as in several respects, $T^*\odot T$ is \bf{not} just the usual tensoring of \nbd{\sB^a(E)}linear maps on internal tensor products of correspondences. Firstly, $T$ is left linear but, in general, not bilinear. (If $T$ was bilinear, it was a \nbd{\tau}isometry.) Secondly, $F^*$ is a Banach right \nbd{\sB^a(E)}module for which $T^*$ is right \nbd{\sB^a(E)}linear, but $F^*$ is not a Hilbert \nbd{\sB^a(E)}module. So, thirdly, $F^*\odot F$ is not an internal tensor product over $\sB^a(E)$.

The proof of boundedness can be done by appealing to the module Haagerup tensor product and Blecher's result \cite[Theorem 4.3]{Ble97} that the internal tensor product  of correspondences is completely isometricly the same as their module Haagerup tensor product. (Indeed, the universal property of the module Haagerup tensor product guarantees that the map $T^*\odot T$ between the module Haagerup tensor norms on the tensor products $E^*\uodot E$ and $F^*\uodot F$ over $\sB^a(E)$ is completely bounded with $\norm{T^*\odot T}_{cb}\le\norm{T^*}_{cb}\norm{T}_{cb}$. The Haagerup seminorm on $F^*\otimes F$ with amalgamation over $\sB^a(E)$, which is homomorphic to a subset of $\sB^a(F)$, is bigger than the Haagerup seminorm with amalgamation over $\sB^a(F)$. So, together with Blecher's result we get that the CB-norm of $\tau$ as map between the internal tensor products is not bigger than $\norm{T^*}_{cb}\norm{T}_{cb}$.) But we prefer to give a direct independent proof.

Let $u=\sum\limits_{i=1}^{n} x_i^*\odot y_i=\sum\limits_{i=1}^{n}\AB{x_i, y_i}\in E^*\uodot E=\ls\AB{E,E}$. For the elements $X^n$ and $Y^n$ in $E^n$ with entries $x_i$ and $y_i$, respectively, this reads $u=\AB{X^n,Y^n}$. We get $(T^*\odot T)(u)=\AB{T^n(X^n),T^n(Y^n)}$. Consequently,
\beqn{
\norm{(T^*\odot T)(u)}
~=~
\norm{\AB{T^n(X^n),T^n(Y^n)}}
~\le~
\norm{T^n}^2\norm{X^n}\norm{Y^n}
~\le~
\norm{T}_{cb}^2\norm{X^n}\norm{Y^n}.
}\eeqn
If, for any $\ve> 0$, we can find $X_\ve$ and $Y_\ve$ in $E^n$ such that $\AB{X_\ve,Y_\ve}=u$ and $\norm{X_\ve}\norm{Y_\ve}\leq \norm{u}+\ve$, then we obtain
\beqn{
\norm{(T^*\odot T)(u)}
~\le~
\norm{T}_{cb}^2\norm{X_\ve}\norm{Y_\ve}
~\le~
\norm{T}_{cb}^2(\norm{u}+\ve),
}\eeqn
and further $\norm{T^*\odot T}\le\norm{T}_{cb}^2$, by letting $\ve\to0$.

For showing that this is possible, we recall the following well-known result. (See, for instance, Lance \cite[Lemma 4.4]{Lan95}.)

\blem\label{lancelem}
For every element $x$ in a Hilbert $\cB$-module $E$ and for every $\alpha\in(0,1)$ there is an element $w_\alpha\in E$ such that $x=w_\alpha\abs{x}^{\alpha}$.
\elem

The proof in \cite{Lan95} shows that $w_\alpha$ can be chosen in the Hilbert \nbd{C^*(\abs{x})}module $\ol{xC^*(\abs{x})}$, which is isomorphic to $C^*(\abs{x})$ via $x\mapsto\abs{x}$. Since $\abs{x}^\alpha$ is strictly positive in the \nbd{C^*}algebra $C^*(\abs{x})$, the element $w_\alpha\in\ol{xC^*(\abs{x})}$ is unique and, obviously, when represented in $C^*(\abs{x})$ it is $w_\alpha=\abs{x}^{1-\alpha}$.

\bcor
Let $E$ be a Hilbert $\cB$-module and let $F$ be a \nbd{\cB}\nbd{\cC}correspondence. Choose $x\in E$, $y\in F$ and put $u:=x\odot y$. Then for every $\ve>0$, there exist $x_\ve\in E$ and $y_\ve\in F$ such that $x_\ve\odot y_\ve=u$ and
\beqn{
\norm{x_\ve}\norm{y_\ve}\leq \norm{u}+\ve,
}\eeqn
that is, $\norm{x\odot y}=\inf\bCB{\,\norm{x'}\norm{y'}\colon x'\in E,y'\in F,x'\odot y'=x\odot y\,}$.
\ecor

\proof
We have $u=x\odot y=w_\alpha\odot\abs{x}^\alpha y$ so that
\beqn{
\norm{u}
~\le~
\norm{w_\alpha}~\norm{\,\abs{x}^\alpha y\,}
~\xrightarrow{~\alpha\to1~}~
1\cdot\norm{\,\abs{x}y\,}
~=~
\norm{x\odot y}
~=~
\norm{u},
}\eeqn
since $\norm{w_\alpha}=\sup_{\lambda\in\SB{0,\norm{x}}}\lambda^{1-\alpha}=\norm{x}^{1-\alpha}\to1$, and since $\abs{x}^\alpha$ converges in norm to $\abs{x}$.\qed

\lf
With the proof of this corollary we did not only \bf{conclude} the proof of \ref{3} $\Rightarrow$ \ref{1}, but also the \bf{proof of Theorem \ref{mainthm}}.\qedsymbol

\bitemp[Corollary {\protect\cite[Theorem 4.3]{Ble97}}.~]
The internal tensor product norm of $u\in E\uodot F$ is
\beq{\label{Haagerupnorm}
\norm{u}
~=~
\inf\BCB{\,\norm{X_n}\norm{Y^n}\colon ~n\in\N,~X_n\in E_n,~Y^n\in F^n,~X_n\odot Y^n=u\,},
}\eeq
with the row space $E_n:=M_{1,n}(E)$ and the internal tensor product $X_n\odot Y^n$ over $M_n(\cB)$. That is, the internal tensor product norm coincides with the module Haagerup tensor product norm (which is defined by \eqref{Haagerupnorm}). Moreover, since $M_n(E\odot F)$ is isomorphic to the internal tensor product $M_n(E)\odot M_n(F)$, the internal tensor product is completely isometricly isomorphic to the module Haagerup tensor product.
\eitemp

\lf
After this digression on the Haagerup tensor product, let us return to maps fulfilling \ref{3}. However, we weaken the conditions a bit. Firstly, we replace $F_T$ with $F$, so that now $F$ is a \nbd{\sB^a(E)}\nbd{\cC}correspondence fulfilling $T(ax)=aT(x)=:\vt(a)T(x)$. We still may define the map $T^*\odot T$ on $E^*\uodot E=\ls\AB{E,E}$, and if $T$ is CB, everything goes as before. Secondly, we wish to weaken the boundedness condition on $T$. We know from Example \ref{counterex} that if $\cB_E$ is nonunital, the CB-condition is indispensable. So, suppose that $E$ is full and that $\cB=\cB_E$ is unital.

\bob\label{uniob}
In the prescribed situation, suppose $E$ has a unit vector $\xi$. In that case, $\tau:=T^*\odot T$ is defined on all $\cB=\AB{\xi,\xi}\cB\subset E^*\uodot E\subset\cB$. Since $\tau(b^*b)=\tau(b^*\AB{\xi,\xi}b)=\AB{(T(\xi b),T(\xi b)}$ is positive, $\tau$ is bounded by $\norm{\tau(\U)}$. From $T(x)=T(x\AB{\xi,\xi})=(x\xi^*)T(\xi)$, we conclude that $\norm{T(x)}^2\le\norm{x}^2\norm{\tau(\U)}$. (This is the same trick in Remark \ref{unirem} that allowed to show that a map $T\colon E\rightarrow F$ fulfilling \ref{3} without boundedness and linearity, is linear provided $E$ has a unit vector $\xi$.)

Even if $E$ has no unit vector but $\cB=\cB_E$ still is unital, then a well-known result asserts that there is a number $n\in\N$ such that $E^n$ has a unit vector, say, $\xi^n$. (See Skeide \cite[Lemma 3.2]{Ske09} for a proof.) If $T$ is linear, then $T^n\colon E^n\rightarrow F^n$ fulfills \ref{3} without boundedness. By the preceding paragraph, $T^n$, and \it{a fortiori} $T$, is bounded by $\sqrt{\norm{\tau(\U)}}$ with the same $\tau$ as obtained from $T$.

Finally, $(T^n)_m=T_{mn,m}\colon M_m(E^n)\rightarrow M_m(F^n)$ is bounded by $\sqrt{\norm{\tau}}$, since $M_m(E^n)$ has a unit vector (with entries $\xi^n$ in the diagonal) and $\norm{\tau_m(\U_m)}=\norm{\tau(\U)}$. So, $T^n$, and \it{a fortiori} $T$, is completely bounded by $\sqrt{\norm{\tau}}$.
\eob

The last missing piece in the proof of Theorem \ref{gentauthm} is the following lemma. We obtain it as a corollary of the proof of Lemma \ref{tauCPlem}.

\blem\label{maxnlem}
$\norm{T}_{cb}\ge\sqrt{\norm{\tau}}$.
\elem

\proof
Let $bb^*$ be in the unitball of $\cB$ such that $\norm{\tau(bb^*)}\approx\norm{\tau}$. By the proof of Lemma \ref{tauCPlem}, there exist $n\in\N$ and $X^n\in E^n$ with $\norm{X^n}\le\norm{b}$ such that $\AB{X^n,X^n}\approx bb^*$ and $\AB{T^n(X^n),T^n(X^n)}\approx\tau(\AB{X^n,X^n})$. So, $\norm{\tau}\approx\norm{\AB{T^n(X^n),T^n(X^n)}}\le\norm{T^n}^2\le\norm{T}_{cb}^2$.\qed

\section{CP-Extendable maps: The KSGNS-con\-struc\-tion revisited}\label{CPrSEC}

In \ref{1} $\Rightarrow$ \ref{2} we have written down the (strict unital) homomorphism $\vt\colon\sB^a(E)\rightarrow\sB^a(F_T)$ in the form $\vt:=v(\bullet\odot\id_\sF)v^*$ with the unitary $v\colon E\odot\sF\rightarrow F_T$ granted by the theorem in \cite{Ske12a}. Then we have shown that the block-wise map $\sT:=\rtMatrix{\tau&T^*\\T&\vt}$ is completely positive, by writing it as $\Xi^*(\bullet\odot\id_\sF)\Xi$ with a \it{diagonal} map $\Xi\in\sB^a\rtMatrix{\tMatrix{\cC\\F_T}\,,\,\tMatrix{\cB\\E}\odot\sF}$. (Recall that it was necessary to unitalize $\tau$ if $\cB$ was nonunital.) We wish to illustrate that these forms for $\vt$ and $\sT$ are not accidental, but they actually are characteristic for all strictly CP-extendable maps $T$.

Let $E$ be a Hilbert \nbd{\cB}module, let $F$ be a Hilbert \nbd{\cC}module, and let $\sT\colon\sB^a(E)\rightarrow\sB^a(F)$ a CP-map. Denote by $(\sE,\Xi)$ the GNS-construction for $\sT$. Like every Hilbert \nbd{\sB^a(F)}module, we may embed $\sE$ into $\sB^a(F,\sE\odot F)$ by identifying $X\in\sE$ with the map $X\odot\id_F\colon y\mapsto X\odot y$ and adjoint $X^*\odot\id_F\colon X'\odot y\mapsto\AB{X,X'}y$. So, $\sT(a)=\Xi^*(a\odot\id_F)\Xi$ where $a\in\sB^a(E)$ acts by the canonical left action on the factor $\sE$ of $\sE\odot F$.

\blem\label{strGNSlem}
The following conditions are equivalent:
\begin{enumerate}
\item\label{str1}
$\sT$ is \hl{strict}, that is, bounded strictly converging nets in $\sB^a(E)$ are sent to strictly converging nets in $\sB^a(F)$.

\item\label{str2}
The action of $\sK(E)$ on the \nbd{\sB^a(E)}\nbd{\cC}correspondence $\sE\odot F$ is nondegenerate.

\item\label{str3}
The left action of the \nbd{\sB^a(E)}\nbd{\cC}correspondence $\sE\odot F$ defines a strict homomorphism.
\end{enumerate}
\elem

\proof
Recall that a correspondence, by definition, has nondegenerate left action, so that $\id_E$ acts as identity. It is well-known (and easy to show) that \ref{str2} and \ref{str3} are equivalent for every \nbd{\sB^a(E)}\nbd{\cC}correspondence. (Indeed, since a bounded approximate unit for $\sK(E)$ converges strictly to $\id_E$, for a strict left action the compacts must act nondegenerately. And if $\sK(E)$ acts nondegenerately, then this action extends to a unique action of all $\sB^a(E)$ that is strict, automatically. See Lance \cite[Proposition 5.8]{Lan95} or the proof of Muhly, Skeide, and Solel \cite[Corollary 1.20]{MSS06}.) Recall, also, that on bounded subsets, strict and \nbd{*}strong topology coincide. (See \cite[Proposition 8.1]{Lan95}.)

Now, if the left action of $\sE\odot F$ is strict, then for every bounded net $\bfam{a_\lambda}_{\lambda\in\Lambda}$ converging strictly to $a$, we have that $(a_\lambda\odot\id_F)(\Xi\odot y)$ converges to $(a\odot\id_F)(\Xi\odot y)$, and likewise for $a_\lambda^*$. In other words, $\Xi^*(a_\lambda\odot\id_F)\Xi$ converges \nbd{*}strongly, hence, strictly to $\Xi^*(a\odot\id_F)\Xi$. So, \ref{3} $\Rightarrow$ \ref{1}.

Conversely, suppose $\sT$ is strict, and choose a bounded approximate unit $\bfam{u_\lambda}_{\lambda\in\Lambda}$ for $\sK(E)$. Then for every element $a\Xi\odot y$ from the total subset $\sB^a(E)\Xi\odot F$ of $\sE\odot F$, we have
\beqn{
\abs{(u_\lambda a-a)\Xi\odot y}^2
~=~
\AB{y,\sT((u_\lambda a-a)^*(u_\lambda a-a))y}
~\longrightarrow~
0,
}\eeqn
so that $\lim_\lambda(u_\lambda\odot\id_F)(a\Xi\odot y)=\lim_\lambda u_\lambda a\Xi\odot y=a\Xi\odot y$. This shows \ref{1} $\Rightarrow$ \ref{2}.\qed

\lf
We now define the \nbd{\cB}\nbd{\cC}correspondence $\sF:=E^*\odot\sE\odot F$. If $\sT$ is strict so that $E\odot E^*\cong\sK(E)$ acts nondegenerately on $\sE\odot F$, then the string
\beqn{
\sE\odot F
~=~
\cls\sK(E)(\sE\odot F)
~\cong~
\sK(E)\odot(\sE\odot F)
~\cong~
(E\odot E^*)\odot(\sE\odot F)
~=~
E\odot(E^*\odot\sE\odot F)
~=~
E\odot\sF
}\eeqn
of (canonical) identifications proves that the map $(x'x^*)(X\odot y)\mapsto x'\odot(x^*\odot X\odot y)$ defines an isomorphism $\sE\odot F\rightarrow E\odot\sF$ of \nbd{\sB^a(E)}\nbd{\cC}correspondences. To obtain the following theorem, we simply have to put the preceding considerations together.%
\footnote{\label{strFN}
This way to construct the \nbd{\cB}\nbd{\cC}correspondence $\sF$ from a \nbd{\sB^a(E)}\nbd{\sB^a(F)}correspondence is, actually, from Bhat, Liebscher, and Skeide \cite[Section 3]{BLS08}. There, however, it is incorrectly claimed that the GNS-correspondence of a strict CP-map has strict left action. (This is false, in general, as the maps $\sT=\id_{\sB^a(E)}$ shows. The results in \cite{BLS08} are, however, correct, as strictness is never used for $\sE$ but always only in the combination as tensor product $\sE\odot F$.) For that reason, we preferred to discuss this here carefully, including also the precise statements in Lemma \ref{strGNSlem}.
}

\bthm\label{strCPthm}
Let $E$ be a Hilbert \nbd{\cB}module, let $F$ be a Hilbert \nbd{\cC}module, and suppose that $\sT\colon\sB^a(E)\rightarrow\sB^a(F)$ is a strict CP-map. Then there exist a \nbd{\cB}\nbd{\cC}correspondence $\sF$ and a map $\Xi\in\sB^a(F,E\odot\sF)$ such that $\Xi^*(\bullet\odot\id_\sF)\Xi=\sT$.
\ethm

\brem
For $E=\cB$ so that $\sB^a(\cB)=M(\cB)$, the \it{multiplier algebra} of $\cB$, this result is known as \hl{KSGNS-construction} for a strict CP-map from $\cB$ into $\sB^a(F)$ (Kasparov \cite{Kas80}); see Lance \cite[Theorem 5.6]{Lan95}. One may consider Theorem \ref{strCPthm} as a consequence of the KSGNS-construction applied to $\sT\upharpoonright\sK(E)$ and the representation theory of $\sB^a(E)$ from Muhly, Skeide, and Solel \cite{MSS06}. Effectively, when $\sT$ is a strict unital homomorphism, so that $\sE={_\sT}\sB^a(F)$ and $\sF:=E^*\odot\sE\odot F=E^*\odot_\sT F$, the theorem (and its proof) specialize to \cite[Theorem 1.4]{MSS06} (and its proof). We like to view Theorem \ref{strCPthm} as a joint generalization of the KSGNS-construction \bf{and} of the representation theory, and the rapid joint proof shows that this point of view is an advantage.
\erem

\bob
Like with all GNS- and Stinespring type constructions, also here we have suitable uniqueness statements. The GNS-correspondence $\sE$ together with the cyclicity condition $\sE=\cls\sB^a(E)\Xi\sB^a(F)$ is unique up to (cyclic-vector-intertwining) isomorphism of correspondences. Of course, this turns over to $\sE\odot F$ with the cyclic map $\Xi\in\sB^a(F,\sE\odot F)$ as with Stinespring construction (as mentioned many times in the sequel of Bhat and Skeide \cite[Example 2.16]{BhSk00} when $F=H$ is a Hilbert space). As for uniqueness of $\sF$, this requires fullness of $E$. Indeed, since $\sE\odot F$ with its action of $\sB^a(E)$ is determined up unitary equivalence, \cite[Theorem 1.8 and Remark 1.9]{MSS06} tell us that $\sF$ is unique if $E$ is full, and and that $\sF$ may fail to be unique if $E$ is not full.
\eob

\bcor
Suppose $E=\rtMatrix{E_1\\E_2}$ and $F=\rtMatrix{F_1\\F_2}$. Then the strict CP-map $\sT$ acts block-wise from $\sB^a(E)=\rtMatrix{\sB^a(E_1)&\sB^a(E_2,E_1)\\\sB^a(E_1,E_2)&\sB^a(E_2)}$ to $\sB^a(F)=\rtMatrix{\sB^a(F_1)&\sB^a(F_2,F_1)\\\sB^a(F_1,F_2)&\sB^a(F_2)}$ if and only if the map $\Xi$ in Theorem \ref{strCPthm} has the diagonal form $\Xi=\rtMatrix{\xi_1~&\\&~\xi_2}$
\ecor

We skip the simple proof.

Now, suppose $\sT=\rtMatrix{\tau&T^*\\T&\vt}\colon\rtMatrix{\cB&E^*\\E&\sB^a(E)}\rightarrow\rtMatrix{\cC&F^*\\F&\sB^a(F)}$ is a block-wise CP-map with strict \nbd{22}corner $\vt$. There is no harm in assuming that $\cC$ is unital. And if $\cB$ is not unital, unitalize $\tau$. For unital $\cB$, the extended linking algebra is $\sB^a\rtMatrix{\cB\\E}$ and the strict topology of all corners but $\sB^a(E)$, coincides with the norm topology. Therefore, $\sT$ is strict. So, except for the possibly necessary unitalization, we see that the form we used in the proof \ref{1} $\Rightarrow$ \ref{2} to establish that the constructed $\sT$ is completely positive, actually, is also necessary. (If unitalization is necessary, then $\xi_1$ is an element of a \nbd{\wt{\cB}}\nbd{\wt{\cC}}correspondence.) We arrive at the factorization theorem for strictly CP-extendable maps, which is the analogue to the theorem in Skeide \cite{Ske12a}.

\bthm\label{CPdecthm}
Let $\cB$ be a unital \nbd{C^*}algebra and let $\cC$ be a \nbd{C^*}algebra. Then for a map $T$ from a Hilbert \nbd{\cB}module $E$ to a Hilbert \nbd{\cC}module $F$ the following conditions are equivalent:
\begin{enumerate}
\item
$T$ admits a strict block-wise extension to a CP-map $\sT=\rtMatrix{\tau&T^*\\T&\vt}\colon\rtMatrix{\cB&E^*\\E&\sB^a(E)}\rightarrow\rtMatrix{\cC&F^*\\F&\sB^a(F)}$.

\item
There exist a \nbd{\cB}\nbd{\cC}correspondence $\sF$, an element $\xi_1\in\sF$ and a map $\xi_2\in\sB^a(F,E\odot\sF)$ such that $T=\xi_2^*(\bullet\odot\xi_1)$.
\end{enumerate}
\ethm
As for a criterion that consists in looking just at $T$, we are reluctant to expect too much. Clearly, such a $T$ must be completely bounded. If $T$ is completely bounded, by appropriate application of Paulsen \cite[Lemma 7.1]{Pau86}, $T$ should extend to the operator system $\rtMatrix{\C\U&E^*\\E&\C\sid_E}\subset\rtMatrix{\cB&E^*\\E&\sB^a(E)}$. But to extend this further, we would have to tackle problems like extending CP-maps from an operator systems to the \nbd{C^*}algebra containing it. We do not know if the special algebraic structure will allow to find a solution to our specific problem. But, in general, existence of such extensions is only granted if the codomain is an injective \nbd{C^*}algebra. We do not follow the question in these notes.

We close this section with an alternative proof of \ref{2} $\Rightarrow$ \ref{1} in Theorem \ref{mainthm}.

\bcor \label{2->1cor}
In the situation of \ref{2} of Theorem \ref{mainthm}, $T$ is a \nbd{T^*\odot T}map.
\ecor

\proof
Recall that the proof \ref{2} $\Rightarrow$ \ref{3} shows us that $\vt$ is unital and strict. Unitalizing if necessary, we get $\xi_1$ and $\xi_2$. Since $\vt$ is a unital homomorphism, $\xi_2$ must be an isometry with $\xi_2\xi_2^*$ commuting with all $a\odot\id_\sF$. This together with $\cls(\sB^a(E)\odot\id_\sF)\xi_2 F_T=E\odot\sF$, implies that $\xi_2$ is unitary. We get $\norm{\AB{T^n(X^n),T^n(X'^n)}}=\norm{\AB{X^n\odot\xi_1,X'^n\odot\xi_1}}\le\norm{\tau}\norm{\AB{X^n,X'^n}}^2$, so $T^*\odot T$ is bounded.\qed

\lf
We think that it is the class of strictly CP-extendable maps that truly merits to be called \hl{CP-maps} between Hilbert modules, and not the more restricted class of CP-H-extendable maps.

\newpage

\section{CPH-semigroups}\label{CPHsgSEC}

In the preceding sections we have seen when a map $T$ from a full Hilbert \nbd{\cB}module $E$ to a Hilbert \nbd{\cC}module $F$ is a \nbd{\tau}map for some CP-map $\tau$ from $\cB$ to $\cC$: If and only if it is CP-H-extendable, that is, if and only if it is a CPH-map into the Hilbert \nbd{\cC}submodule generated by $T(E)$, $F_T$. If $E$ is not full, then this may be repaired easily by making $\cB$ smaller. If a CP-H-extendable map is not a CPH-map, then this may be repaired easily by making $F$ smaller. In fact, we have seen that replacing $F$ with $F_T$, we turn $T$ even into a strictly CPH$_0$-map $E\rightarrow F_T$. In that case, the CPH-extension $\sT=\rtMatrix{\tau&T^*\\T&\vt}$ is even unique.

Similarly, the conditions in Theorem \ref{mainthm} tell when a semigroup $T=\bfam{T_t}_{t\in\R_+}$ of maps $T_t$ on a full Hilbert \nbd{\cB}module $E$ is \hl{CP-H-extendable}, that is, when each map $T_t$ is CP-H-extendable. In this case, it is even clear that the (unique) maps $\tau_t$ turning the $T_t$ into \nbd{\tau_t}maps, form a CP-semigroup $\tau$ on $\cB$. However, the situation is considerably different, when we ask if the $T_t$ are actually CPH-maps. In the sequel, we shall see that no such semigroup will ever fulfill $E=E_{T_t}$ for all $t$, unless all $\tau_t$ are homomorphic (see Observation \ref{notETob}) and, therefore, the $T_t$ are ternary homomorphisms. We shall see that we may replace the unfulfillable condition $E=E_{T_t}$ with a weaker minimality condition (Definition \ref{minCPHdefi}) involving the whole semigroup, which also will guarantee existence of (unique) strictly CPH$_0$-extensions $\sT_t=\rtMatrix{\tau_t&T_t^*\\T_t&\vt_t}$ which even form a semigroup themselves. Understanding this, requires results from Bhat and Skeide \cite{BhSk00} about the GNS-product system of a CP-semigroup (replacing Paschke's GNS-construction for a single CP-map) and about the relation between product systems and strict \nbd{E_0}semigroups on $\sB^a(E)$ from Skeide \cite{Ske02,Ske09}. The construction of minimal CPH-semigroups involves results about existence of \nbd{E_0}semigroups for product systems from Skeide \cite{Ske06,Ske07,Ske08p1}.

\lf
Let us first fix the sort of semigroup we wish to look at. Recall that an \hl{\nbd{E}semigroup} is a semigroup of endomorphisms on a \nbd{*}algebra, and that an \hl{\nbd{E_0}semigroup} is a semigroup of unital endomorphisms on a unital \nbd{*}algebra.

\bdefi\label{CP(H)sgdef}
A semigroup $T=\bfam{T_t}_{t\in\R_+}$ of maps $T_t\colon E\rightarrow E$ on a Hilbert \nbd{\cB}module $E$ is
\begin{enumerate}
\item\label{sg1}
a (\hl{strictly}) \hl{CP-semigroup} on $E$ if it extends to a CP-semigroup $\sT=\bfam{\sT_t}_{t\in\R_+}$ of maps $\sT_t=\rtMatrix{\tau_t&T_t^*\\T_t&\vt_t}$ acting block-wise on the extended linking algebra of $E$ (with strict $\vt_t$);

\item\label{sg2}
a (\hl{strictly}) \hl{CPH$(_0)$-semigroup} on $E$ if it is a (strictly) CP-semigroup where the $\vt_t$ can be chosen to form an \nbd{E(_0)}semigroup and where the $\tau_t$ can be chosen such that each $T_t$ is a \nbd{\tau_t}map.
\end{enumerate}
\edefi

\bob\label{unistriob}
In the sequel, frequently the results will depend on that $\cB$ is a unital \nbd{C^*}al\-gebra. Recall that, by the discussion preceding Theorem \ref{CPdecthm}, in this case $T$ being a strictly CP-semigroup (and so forth) on a Hilbert \nbd{\cB}module, simply means that each $\sT_t$ is strict. In that case, we will just say, $T$ is a strict CP-semigroup (and so forth).
\eob

In the sequel, we shall address the following problems: We give a version of the decomposition in Theorem \ref{CPdecthm} for strict CP-semigroups; Theorem \ref{CPsgdthm}. In order to prepare better for the case of CPH-semigroups, we are forced to be more specific than in Section \ref{CPrSEC}; see the extensive Observation \ref{PSsTob}. Then, we examine to what extent this version for CPH-semigroups corresponds to the single map results from Skeide \cite{Ske12a} and Theorem \ref{mainthm}. The version in Theorem \ref{sgSkethm} for CP-H-extendable semigroups of the single map result in \cite{Ske12a} is preliminary for the result Theorem \ref{sgthm} on CPH-semigroups. The latter result parallels rather Theorem \ref{CPsgdthm} (hypothesizing that there is CPH-extension of the CP-H-extendable semigroup $T$), than proving existence of a CPH-extension, as in Theorem \ref{mainthm}, from CP-H-extendability under (here, unfulfillable) cyclicity conditions. The results that parallels Theorem \ref{mainthm} most, is Theorem \ref{minCPH0thm} on \it{minimal} CP-H-extendable semigroups on full Hilbert modules over unital \nbd{C^*}algebras. The minimality condition in \eqref{Tmin} limits this theorem automatically to the case where the associated CP-semigroups have full \it{GNS-systems}. In this case, however, we can prove existence (based on the corresponding existence results of \nbd{E_0}semigroups for such full product systems). We show that all minimal CP-H-extendable semigroups on a fixed full Hilbert \nbd{\cB}module and \it{associated} with a fixed CP-semigroup on $\cB$, are \it{cocycle equivalent}.

\lf
We start by discussing what we can say about strict CP-semigroups on $\sB^a(E)$ in general. As the basis for Theorem \ref{CPdecthm} and the other results in Section \ref{CPrSEC} is Paschke's GNS-construction for a single CP-map $\tau$, here we will need the version of the GNS-construction for CP-semigroups from Bhat and Skeide \cite{BhSk00}.

Let $\tau=\bfam{\tau_t}_{t\in\R_+}$ be a CP-semigroup on a unital \nbd{C^*}algebra $\cB$. Bhat and Skeide \cite[Section 4]{BhSk00} provide the following:
\begin{itemize}
\item
A \hl{product system} $E^\odot=\bfam{E_t}_{t\in\R_+}$ of \nbd{\cB}correspondences. That is, $E_0=\cB$ (the trivial \nbd{\cB}correspondence), and there are bilinear unitaries $u_{s,t}\colon E_s\odot E_t\rightarrow E_{s+t}$ such that the \hl{product} $(x_s,y_t)\mapsto x_sy_t:=u_{s,t}(x_s\odot y_t)$ is associative and such that $u_{0,t}$ and $u_{t,0}$ are left and right action, respectively, of $\cB=E_0$ on $E_t$.

\item
A \hl{unit} $\xi^\odot=\bfam{\xi_t}_{t\in\R_+}$ (that is, the elements $\xi_t\in E_t$ fulfill $\xi_0=\U$ and $\xi_s\xi_t=\xi_{s+t}$), such that
\beqn{
\tau_t
~=~
\AB{\xi_t,\bullet\xi_t},
}\eeqn
and the smallest product subsystem of $E^\odot$ containing $\xi^\odot$ is $E^\odot$. (The pair $(E^\odot,\xi^\odot)$ is determined by these properties up to unit-preserving isomorphism, and we refer to it as the \hl{GNS-construction} for $\tau$ with \hl{GNS-system} $E^\odot$ and \hl{cyclic unit} $\xi^\odot$.)

\item
If $E^\odot$ is not minimal, then the subcorrespondences
\beq{\label{GNSform}
E_t
~=~
\cls\bCB{b_n\xi_{t_n}\ldots b_1\xi_{t_1}b_0\colon t_n+\ldots+t_1=t}.
}\eeq
of $E_t$ form a product subsystem of $E^\odot$ that is isomorphic to the GNS-system.
\end{itemize}
Now let $\sT=\bfam{\sT_t}_{t\in\R_+}$ be a CP-semigroup on the unital \nbd{C^*}algebra $\sB^a(E)$. (For this, $\cB$ need not even be unital.) Denote by $\sE^\odot=\bfam{\sE_t}_{t\in\R_+}$ its GNS-system and by $\Xi^\odot=\bfam{\Xi_t}_{t\in\R_+}$ its cyclic unit. Like in Lemma \ref{strGNSlem}, the semigroup $\sT$ is strict if and only if the correspondences have strict left action. (To see this, it is crucial to know the form in \eqref{GNSform} of a typical element in the GNS-system.) Like in Theorem \ref{strCPthm}, if $\sT$ is strict, we get \nbd{\cB}correspondences $E_t:=E^*\odot\sE_t\odot E$ (actually, \nbd{\cB_E}correspondences if $E$ is not full) and maps in $\sB^a(E,E\odot E_t)$, also denoted by $\Xi_t$, such that $\sT_t=\Xi_t^*(\bullet\odot\id_t)\Xi_t$.

In addition to the properties discussed in Section \ref{CPrSEC}, we see that the $E_t$ form a product system of \nbd{\cB_E}correspondences via
\beqn{
u_{s,t}
\colon
E_s\odot E_t
~=~
E^*\odot\sE_s\odot E\odot E^*\odot\sE_t\odot E
~\rightarrow~
E^*\odot\sE_s\odot\sE_t\odot E
~\rightarrow~
E^*\odot\sE_{s+t}\odot E
~=~
E_{s+t},
}\eeqn
and the $\Xi_t$ compose as
\beq{\label{Xiunit}
\Xi_{s+t}
~=~
(\id_E\odot u_{s,t})(\Xi_s\odot\id_t)\Xi_t.
}\eeq
(Note that, modulo the flaw in Bhat, Liebscher, and Skeide \cite{BLS08} regarding strictness of the GNS-construction mentioned in Footnote \ref{strFN}, all this has already been discussed in \cite{Ske09} and in \cite{BLS08}.) Of course, every product system $E^\odot$ with a family of maps $\Xi_t\in\sB^a(E,E\odot E_t)$ satisfying \eqref{Xiunit}, defines a strict CP-semigroup $\sT$ on $\sB^a(E)$ by setting $\sT_t:=\Xi_t^*(\bullet\odot\id_t)\Xi_t$. But only if $E^\odot$ and the $\Xi_t$ arise as described, we will speak of the \hl{product system} of $\sT$.

It is worth to collect some properties of the product system $E^\odot$ of $\sT$ and the $\Xi_t$.

\bob\label{PSsTob}
\begin{enumerate}
\item\label{PSsT1}
Recall, that $a\odot\id_t\in\sB^a(E\odot E_t)$, when composed with an element $X_t\in\sE_t\subset\sB^a(E,E\odot E_t)$, is nothing but the left action of $a\in\sB^a(E)$ on $X_t\in\sE_t$. Therefore, it is sometimes convenient to write $aX_t$ instead of $(a\odot\id_t)X_t$.
Note, too, that by the way how $\sE_t\odot E$ is canonically identified with $E\odot E_t=E\odot E^*\odot\sE_t\odot E=\cls\sK(E)\sE_t\odot E$, we get
\beqn{
x\odot(y^*\odot X_t\odot z)
~=~
(xy^*)X_tz
~\in~
E\odot E_t.
}\eeqn

\item\label{PSsT2}
By the way how $\sE_t$ is generated from $\Xi^\odot$ as expressed in \eqref{GNSform}, it follows from \eqref{Xiunit} that
\beqn{
E\odot E_t
~=~
\cls\BCB{(a_n\odot\id_t)(\Xi_{t_n}a_{n-1}\odot u_{t_{n-1},t_{n-2}+\ldots+t_1})\ldots(\Xi_{t_3}a_2\odot u_{t_2,t_1})(\Xi_{t_2}a_1\odot\id_{t_1})\Xi_{t_1}x}.
}\eeqn
(If $E$ is full, one may show that $E^\odot$ and the $\Xi_t$ are determined uniquely by $\sT$ and that cyclicity condition. But we do not address uniqueness here.) Observe that it suffices to choose the $a_k$, which \it{a priori} run over all $\sB^a(E)$, only from the rank-one operators. Doing so and tensoring with $E^*$ from the left, we get
\baln{
E_t
&
~=~
\cls\BCB{(y_n^*\odot\Xi_{t_n}\odot z_n)\ldots(y_1^*\odot\Xi_{t_1}\odot z_1)}
\\
&
~=~
\cls\BCB{((y_n^*\odot\id_{t_n})\Xi_{t_n}z_n)\ldots((y_1^*\odot\id_{t_1})\Xi_{t_1}z_1)}.
}\ealn
This means, $E^\odot$ is generated as a product system by the family of subsets $E^*\odot\Xi_t\odot E$ of $E_t$.

\item\label{PSsT3}
For both exploiting the preceding cyclicity condition and making notationally the connection with the construction of a product system for strict \nbd{E}semigroups on $\sB^a(E)$, it is convenient to replace the maps $\Xi_t$ with their adjoints $v_t:=\Xi_t^*\colon E\odot E_t\rightarrow E$. Using the same product notation $xy_t:=v_t(x\odot y_t)$ as for the $u_{s,t}$, Equation \eqref{Xiunit} transform into the associativity condition
\beq{\label{assoc}
(xy_s)z_t
~=~
x(y_sz_t).
}\eeq
It follows from the cyclicity condition that we know $v_t$ fulfilling associativity if we know each $v_t$ on $E^*\odot\Xi_t E$. In particular, for checking if $v_t$ is an isometry, it suffices to check that each $v_t$ is an isometry on the subset $E^*\odot\Xi_t\odot E$ of $E_t$.

\item\label{PSsT4}
Isometric $v_t:=\Xi_t^*$ arise from \nbd{E}semigroups in the following way. Observe that $\sT$ is a unital if and only if the $v_t$ are coisometries. (If $E$ is full, this means that $E^\odot$ is necessarily full.) $\sT$ is an \nbd{E}semigroup if and only if for each $t$, $v_t^*v_t$ is a projection in the relative commutant of $\sB^a(E)\odot\id_t$ in $\sB^a(E\odot E_t)$. (This happens, for instance, if the $v_t$ are isometries, so that $v_t^*v_t=\id_{E\odot\sid_t}$ commutes with everything.)

Now, if $\sT$ is an \nbd{E}semigroup, then necessarily $(a\odot\id_t)\Xi_t=\Xi_t\sT_t(a)$. Then, in the cyclicity condition all $a_n$ can be put through to the right, where they are applied to $x$, and the remaining $\Xi_{t_k}$, following \eqref{Xiunit}, multiply together to give $\Xi_t$. We conclude that $E\odot E_t=\Xi_tE$. (Since $\Xi_t$ is a partial isometry, no closure of the image of $\Xi_t$ is necessary.) In other words, the $\Xi_t$ are coisometries, that is, the $v_t$ are isometries. So, if $\sT$ is an \nbd{E_0}semigroup, then the $v_t$ are even unitaries.

\item\label{PSsT5}
Since \cite{Ske02} (\cite{Ske09} (preprint 2004) for full $E$ over general \nbd{C^*}algebras), it is known that every strict \nbd{E_0}semigroup (\nbd{E}semigroup) $\sT$ on $\sB^a(E)$ comes along with a product system $E^\odot$ and a family $v_t\colon E\odot E_t\rightarrow E$ of unitaries (adjointable isometries) fulfilling \eqref{assoc}, such that $\sT_t=v_t(\bullet\odot\id_t)v_t^*$. If $E$ is full, since \cite{Ske06,Ske07} this is referred to as \hl{left dilation} (\hl{left semi-dilation} of $E^\odot$ to $E$. It is known that product system and left (semi\nobreakdash-)dilation are essentially unique. (It is part of the extensive \cite[Proposition 6.3]{Ske08p1} to explain in which sense these objects are unique.) In fact, it is not difficult to verify that the left (semi-)dilation constructed above, coincides with the one constructed in \cite{Ske09}. But we do not need this information. If $E$ is not full, then we speak of \hl{left quasi-dilation} (\hl{left quasi-semidilation}).
\end{enumerate}
\eob

\noindent
This lengthy observation prepares the ground for Theorem \ref{sgthm}. But logically it belongs here, where strict CP-semigroups on $\sB^a(E)$ are discussed. Of course, the statements regarding single CP-maps acting block-wise on $\sB^a\rtMatrix{E_1\\E_2}$ remain true for strict CP-semigroups. That is $\Xi_t=\rtMatrix{\xi^1&\\&\xi^2_t}\in\sB^a\rtMatrix{\rtMatrix{E_1\\E_2},\rtMatrix{E_1\\E_2}\odot E_t}$. Again, when $E_1=\cB\ni\U$ and $E_2=E$ (so that $\sB^a\rtMatrix{E_1\\E_2}$ is the extended linking algebra of $E$), we identify $\xi^1_t$ with the elements $\xi_t:=\xi^1_t\U\in E_t$. We put $v_t={\xi^2_t}^*$, getting:

\bthm \label{CPsgdthm}
Let $\cB$ be a unital \nbd{C^*}algebra. Then for a semigroup $T=\bfam{T_t}_{t\in\R_+}$ of maps on a Hilbert \nbd{\cB}module $E$ the following conditions are equivalent:
\begin{enumerate}
\item\label{CPs1}
$T$ is a strict CP-semigroup.

\item\label{CPs2}
There exist a product system $E^\odot=\bfam{E_t}_{t\in\R_+}$ of \nbd{\cB}correspondences, a unit $\xi^\odot$ for $E^\odot$, and a family $\bfam{v_t}_{t\in\R_+}$ of maps $v_t\in\sB^a(E\odot E_t,E)$ fulfilling \eqref{assoc}, such that $T_t=v_t(\bullet\odot\xi_t)$.
\end{enumerate}
\ethm

\brem\label{non1strirem}
Note that $E$ is not required full. But, unitality of $\cB$ enters in two ways. Firstly, $\cB$ must be unital in order to obtain the unit $\xi^\odot$ in the product system. (Recall that the term \it{unit} is not defined if $\cB$ is nonunital.) Secondly and more importantly, the construction of the product system $E^\odot$ starts from a strict CP-semigroup $\sT$ on $\sB^a\rtMatrix{\cB\\E}$. For both facts the fact that $\sB^a\rtMatrix{\cB\\E}$ is the extended linking algebra appearing in the definition of CP-semigroup on $E$ and the fact that the CP-extension $\sT$ of $T$ be strict, it is indispensable that $\cB$ is unital; see Observation \ref{unistriob}.

If $\cB$ is nonunital (for instance, because we wish to consider $E$ as full), then the construction of the product system my be saved provided $T$ really may be extended to a CP-semigroup acting strictly on the bigger algebra $\sB^a\rtMatrix{\cB\\E}$. This requires that the semigroup $T$ itself extends to a semigroup of strict maps on $\sB^a(\cB,E)\supset E$. It also requires that the there there is a strict CP-semigroup on $\sB^a(\cB)=M(\cB)$ extending $\tau$. If all this is fulfilled, then instead of a unit for the product system $E^\odot$ we obtain a family of maps $\xi^1_t\in\sB^a(\cB,E_t)$ fulfilling a condition similar to \eqref{Xiunit}. While a product system can be obtained from $\tau$ on nonunital $\cB$ also when $\tau$ is not strict, existence of the maps $\xi^1_t$ is unresolvably intervowen with strictness of $\tau$.

We do not address these questions here. We just mention that there have already been several instances where such \it{multiplier spaces} $\sB^a(\cB,E_t)$ and their \it{strict} tensor products like $\sB^a(\cB,E_s)\odot^{str}\sB^a(\cB,E_t):=\cls^{str}(E_s\odot\id_t)E_s=\sB^a(\cB,E_s\odot E_t)$ popped up. It would be interesting to formulate a theory of product systems for them, extending what has been said in Skeide \cite[Section 7]{Ske09}. Families of maps like $\Xi_t$ fulfilling \eqref{Xiunit} (and, of course, $\Xi_0=\U\in M(\cB)$) would generalize the concept of unit.
\erem

\bob\label{subob}
Note that $E^\odot$ need not be the GNS-system of $\tau$. Of course, it contains the GNS-system, because it contains the unit $\xi^\odot$ that gives back $\tau$ as $\tau_t=\AB{\xi_t,\bullet\xi_t}$. Also the product system of \nbd{\cB_E}correspondences constructed as before from the strict CP-semigroup $\vt$ on $\sB^a(E)$ given by $\vt_t=v_t(\bullet\odot\id_t)v_t^*$ on $\sB^a(E)$ sits inside $E^\odot$. More precisely, by Observation \ref{PSsTob}\eqref{PSsT2}, $E_t$ is generated by $\rtMatrix{\cB\\E}^*\odot\Xi_t\odot\rtMatrix{\cB\\E}$ and, by diagonality of $\Xi_t$, we have
\beqn{
\rtMatrix{\cB\\E}^*\odot\Xi_t\odot\rtMatrix{\cB\\E}
~=~
\rtMatrix{\cB\\0}^*\odot\rtMatrix{\xi^1_t&0\\0&0}\odot\rtMatrix{\cB\\0}+\rtMatrix{0\\E}^*\odot\rtMatrix{0&0\\0&\xi^2_t}\odot\rtMatrix{0\\E}.
}\eeqn
And $\rtMatrix{c\\0}^*\odot\rtMatrix{\xi^1_t&0\\0&0}\odot\rtMatrix{d\\0}$ can be identified with the element $c^*\xi_td$ from the subset generating the GNS-system of $\tau$, while $\rtMatrix{0\\y}^*\odot\rtMatrix{0&0\\0&\xi^2_t}\odot\rtMatrix{0\\z}$ can be identified with the element $y^*\odot\xi^2_t\odot z$ from the subset generating the product system of $\vt$. It is clear that the product system of $\vt$, consisting of \nbd{\cB_E}correspondences (that may also be viewed as \nbd{\cB}correspondences), must be smaller than $E^\odot$ if $E$ is non-full. But it may be smaller even if $E$ is full. (Think of $E=\cB$, where $\tau$ is the identity and $\vt$ a CP-semigroup with nonfull GNS-system.)
\eob

The situation in this observation, namely, that neither of the product systems of the diagonal corners $\tau$ and $\vt$ need coincide with the product system $E^\odot$ of $\sT$, creates not little discomfort. This improves if $T$ is a strict CPH-semigroup, to which we now gradually switch our attention.

For instance, we know that $\vt=v_t(\bullet\odot\id_t)v_t^*$ is an \nbd{E}semigroup if and only if $v_t^*v_t$ is a projection commuting with $\sB^a(E)\odot\id_t$. We also know that, if $E^\odot$ actually is the product system of $\vt$, then the $v_t$ will be isometries. But, if $E^\odot$ is too big, then there is no \it{a priori} reason, why the $v_t$ should be isometries.

It is one of the scopes of the following theorem to contribute an essential part in the proof that the $v_t$ actually \it{are} isometries. A second scope is to present the semigroup version of the theorem in Skeide \cite{Ske12a}. (This will allow to show that the condition $E=E_{T_t}$ can be rarely fulfilled, and by what it has to be replaced. It will also lead to a notion of \it{minimal} CPH-semigroups.)

\bthm\label{sgSkethm}
Let $\cB$ be a unital \nbd{C^*}algebra and let $T=\bfam{T_t}_{t\in\R_+}$ be a family of maps on a Hilbert \nbd{\cB}module $E$. Then the following conditions are equivalent:
\begin{enumerate}
\item\label{CPHex1}
$T$ is a CP-H-extendable semigroup.

\item\label{CPHex2}
There are a product system $E^\odot$ of \nbd{\cB}correspondences, a unit $\xi^\odot$ for $E^\odot$, and a family of (not necessarily adjointable) isometries $v_t\colon E\odot E_t\rightarrow E$ fulfilling \eqref{assoc}, such that
\beqn{
T_t
~=~
v_t(\bullet\odot\xi_t).
}\eeqn
\end{enumerate}
\ethm

\proof
Of course, given the ingredients in \ref{CPHex2}, the maps $T_t$ defined there are CP-H-extendable. The semigroup property follows from the unit property of $\xi^\odot$ and from \eqref{assoc}. This shows \ref{CPHex1}.

Now let $T$ be a CP-H-extendable semigroup. Denote by $\tau$ a CP-semigroup on $\cB$ such that each $T_t$ is a \nbd{\tau_t}map. Do the GNS-construction for $\tau$ to obtain $(E^\odot,\xi^\odot)$. Recall that $E_t$ is spanned by elements as in \eqref{GNSform}. Then
\baln{
x\odot(b_n\xi_{t_n}\ldots b_1\xi_{t_1})
&
~\longmapsto~
T_{t_1}(T_{t_2}(\ldots T_{t_{n-1}}(T_{t_n}(xb_n)b_{n-1})\ldots b_2)b_1)
&
(t_n+\ldots+t_1
&
~=~
t)
}\ealn
extend to well-defined isometries $v_t\colon E\odot E_t\rightarrow E$ fulfilling all the requirements of \ref{CPHex2}. (In order to compute inner products of two elements of the form as in \eqref{GNSform}, one first has to asure, by splitting pieces $\xi_r$ of the unit suitably into $\xi_{r'}\xi_{r''}$, that both elements belong to the same tuple $t_n+\ldots+t_1=t$. We leave the remaining statements to the reader.)\qed

\lf
Note that $E$ is not required full. (It should be specified that also in this case, by a CP-H-extendable map $T$ on $E$ we mean that $T$ is a CPH-map into $E_T$. Likewise, in the semigroup version it is required that the $\tau_t$ turning $T_t$ into \nbd{\tau_t}maps, form a semigroup.) This is, why $\tau$ is not unique. If we wish to emphasize a fixed CP-semigroup $\tau$, we say $T$ is a CP-H-extendable semigroup \hl{associated} with $\tau$.

The proof also shows that $(E^\odot,\xi^\odot)$ may be chosen to be the GNS-construction for $\tau$. But for the backward direction, this is not necessary.

If, by any chance, we find $E^\odot$ and $\xi^\odot$ such that the $v_t$ can be chosen adjointable (so that they form a left quasi-semidilation), then we get that $T$ is a strict CPH-semigroup. (Define the members of the semigroup $\sT$ in the very same way as the single map $\sT$ in the proof of \ref{1} $\Rightarrow$ \ref{2} in Theorem \ref{mainthm}.) If the $v_t$ can even be chosen unitary (so that they form a left quasi-dilation), then $T$ turns out to be a strict CPH$_0$-semigroup. After Observations \ref{PSsTob} and \ref{subob} and after Theorem \ref{sgSkethm}, we now are prepared to prove the opposite direction, too:

\bthm \label{sgthm}
Let $\cB$ be a unital \nbd{C^*}algebra and let $T=\bfam{T_t}_{t\in\R_+}$ be a family of maps on a Hilbert \nbd{\cB}module $E$. Then the following conditions are equivalent:
\begin{enumerate}
\item\label{sgt1}
$T$ is a strict CPH-semigroup (CPH$_0$-semigroup).

\item\label{sgt2}
There exist a product system $E^\odot$, a unit $\xi^\odot$ for $E^\odot$, and a left quasi-semidilation (a left quasi-dilation) $\bfam{v_t}_{t\in\R_+}$ of $E^\odot$ to $E$, such that $T_t=v_t(\bullet\odot\xi_t)$.
\end{enumerate}
\ethm

\proof
Only the direction \ref{sgt1} $\Rightarrow$ \ref{sgt2} is yet missing. So, let $T$ be a strict CPH-semigroup on $E$ and $\sT$ a suitable strict CPH-extension to the extended linking algebra $\sB^a\rtMatrix{\cB\\E}$ of $E$ and construct everything as for \ref{CPs1} $\Rightarrow$ \ref{CPs2} of Theorem \ref{CPsgdthm}. So, $E_t=\rtMatrix{\cB\\E}^*\odot\sE_t\odot\rtMatrix{\cB\\E}$ and $v_t$ is the (co)restriction of $\Xi_t^*\colon\rtMatrix{\cB\\E}\odot E_t\rightarrow\rtMatrix{\cB\\E}$ to the map ${\xi^2_t}^*\colon E\odot E_t\rightarrow E$. By Observation \ref{subob}, $E_t$ is generated by its subset
\beqn{
\rtMatrix{\cB\\0}^*\!\!\odot\rtMatrix{\xi^1_t&0\\0&0}\odot\rtMatrix{\cB\\0}+\rtMatrix{0\\E}^*\!\!\odot\rtMatrix{0&0\\0&\xi^2_t}\odot\rtMatrix{0\\E},
}\eeqn
and by Observations \ref{PSsTob}\eqref{PSsT4}, it is sufficient to check isometry of $v_t\colon E\odot E_t\rightarrow E$ for $x\odot y_t$ where $y_t$ are chosen from that subset. So, we have to check
\beqn{
\AB{xy_t,x'y'_t}
~=~
\AB{x\odot y_t,x'\odot y'_t}
}\eeqn
where $x,x'\in E$ and $y_t,y'_t\in\rtMatrix{\cB\\0}^*\!\!\odot\rtMatrix{\xi^1_t&0\\0&0}\odot\rtMatrix{\cB\\0}\cup\rtMatrix{0\\E}^*\!\!\odot\rtMatrix{0&0\\0&\xi^2_t}\odot\rtMatrix{0\\E}$. Now, for elements $y_t$ and $y'_t$ in the first set (which generates the GNS-system of $\tau$), Theorem \ref{sgSkethm} tells us that $v_t$ in this case is isometric. For elements $y_t$ and $y'_t$ in the second set (which generates the product system of $\vt$), it is easy to see that the $v_t$ in this case give back the $v_t$ of $\vt$, which, we know, are isometric. So, it remains to check the case where $y_t=\rtMatrix{c\\0}^*\odot\Xi_t\odot\rtMatrix{d\\0}$ is from the first set and $y'_t=\rtMatrix{0\\y}^*\odot\Xi_t\odot\rtMatrix{0\\z}$ is from the second set.

We use all the notation from Observations \ref{PSsTob}\eqref{PSsT1}. Additionally, note that by the proof of Lemma \ref{bilcp}, it follows that
\beqn{
\rtMatrix{0&0\\0&a}\Xi_t
~=~
\rtMatrix{0&0\\0&\sid_E}\Xi_t\rtMatrix{0&0\\0&\vt_t(a)},
~=~
\rtMatrix{0&0\\0&\xi^2_t\vt_t(a)}
}\eeqn
and further
\beqn{
\Xi_t\Xi_t^*\rtMatrix{0&0\\0&a}\Xi_t
~=~
\rtMatrix{0&0\\0&\xi^2_t{\xi^2_t}^*\xi^2_t\vt_t(a)}
~=~
\rtMatrix{0&0\\0&\xi^2_t\vt_t(a)}
~=~
\rtMatrix{0&0\\0&a}\Xi_t,
}\eeqn
because $\xi^2_t$ is a partial isometry. We find
\bmun{
\AB{xy_t,x'y'_t}
~=~
\BAB{\Xi_t^*\family{\rtMatrix{0\\x}\odot\rtMatrix{c\\0}^*\!\!\odot\Xi_t\odot\rtMatrix{d\\0}},\Xi_t^*\family{\rtMatrix{0\\x}\odot\rtMatrix{0\\y}^*\!\!\odot\Xi_t\odot\rtMatrix{0\\z}}}
\\[1ex]
~=~
\BAB{\rtMatrix{0\\x}\odot\rtMatrix{c\\0}^*\!\!\odot\Xi_t\odot\rtMatrix{d\\0},\Xi_t\Xi_t^*\rtMatrix{0&0\\0&xy^*}\Xi_t\rtMatrix{0\\z}}
~=~
\BAB{\rtMatrix{0\\x}\odot\rtMatrix{c\\0}^*\!\!\odot\Xi_t\odot\rtMatrix{d\\0},\rtMatrix{0&0\\0&xy^*}\Xi_t\rtMatrix{0\\z}}
\\[1ex]
~=~
\BAB{\rtMatrix{0\\x}\odot\rtMatrix{c\\0}^*\!\!\odot\Xi_t\odot\rtMatrix{d\\0},\rtMatrix{0\\x}\odot\rtMatrix{0\\y}^*\!\!\odot\Xi_t\odot\rtMatrix{0\\z}}
~=~
\AB{x\odot y_t,x'\odot y'_t},
}\emun
so the $v_t$ are, indeed, isometries. And, of course, $\vt$ is an \nbd{E_0}semigroup if and only if the $v_t$ are unitary.\qed

\lf
Every product system $E^\odot$ can be recovered easily from a strict \nbd{E}semigroup $\vt$ acting on a suitable $E$. Indeed, take $E=L^2(E^\odot)$, the direct integral $\int_0^\infty E_\alpha\,d\alpha$ over the product system. (If $E^\odot$ is just a product system, then we have to stick to the counting measure on $\R_+$, that is, $E=\bigoplus_{t\in\R_+}E_t$. If $E^\odot$ is a \it{continuous product system} in the sense of Skeide \cite[Section 7]{Ske03b}, we take the Lebesgue measure and $E$ is the norm closure of the \it{continuous sections} with compact support.) Then the obvious isomorphism from $E\odot E_t$ onto the submodule $\int_t^\infty E_\alpha\,d\alpha$ of $E$ defines a left semidilation $v_t$, and $\vt$ defined by $\vt_t:=v_t(\bullet\odot\id_t)v_t^*$ has product system $E^\odot$. (Thanks to $E_0=\cB$, the module $E$ is full. For the direct sum this is clear. For the continuous case, fullness follows from fullness of $E_0$ and from existence for every $x_0\in E_0$ of a continuous section assuming that value $x_0$ at $\alpha=0$.)  It is easy to see that for a continuous product system, $\vt$ is strongly continuous. Also, by Skeide \cite[Appendix A.1 ]{Ske08p1} applied to the unitalization of $\tau$, the GNS-system of a strongly continuous and contractive CP-semigroup $\tau$ on a unital \nbd{C^*}algebra is continuous.

The backward implication of Theorem \ref{sgthm} gives the following:

\bcor\label{sgcor}
Let $\tau$ be a (strongly continuous) CP-semigroup (of contractions) on the unital \nbd{C^*}algebra $\cB$. Then there exists a (strongly continuous) CPH-semigroup $T$ on a full Hilbert \nbd{\cB}module associated with $\tau$.
\ecor

If we can construct for the (continuous) GNS-system or any (continuous ) product system containing it an \nbd{E_0}semigroup, then we even get a (strongly continuous) CPH$_0$-semigroup $T$. For the existence results of such \nbd{E_0}semigroups, however, it is indispensable that this product system $E^\odot$ is full. Continuous product systems ($\cB$ unital!) are full; see \cite[Lemma 3.2]{Ske07}. (Note that this is a result that does not hold for von Neumann correspondences.) GNS-Systems of so-called \it{spatial} CP-semigroups (continuous or not) embed into a full product system; see Bhat, Liebscher, and Skeide \cite{BLS10}. (It is strongly full in the von Neumann case; see \cite[Theorem A.15]{Ske08p1}.) Of course, the GNS-system of a Markov semigroup is full. For Markov semigroups, there is an easy way to construct \nbd{E_0}semigroup, to which we will come back in Section \ref{CPHdilSEC}. For nonunital semigroups, we have to stick to the existence result in Skeide \cite{Ske07}, which generalizes to modules the proof in Skeide \cite{Ske06} of Arveson's fundamental results \cite{Arv90} that every product system of Hilbert spaces comes from an \nbd{E_0}semigroup on $\sB(H)$. (The von Neumann case is dealt with in \cite{Ske08p1}.) We see that all Markov semigroups and most CP-semigroups have CPH-semigroups with which they are associated.

\lf
We conclude this section by drawing some consequences from Theorem \ref{sgSkethm}. In particular, we wish to find information how to make sure that a CP-H-extendable semigroup either is a strict CPH-semigroup.

Well, given a unital \nbd{C^*}algebra $\cB$ and a CP-H-extendable semigroup on a Hilbert \nbd{\cB}module $E$ associated with a CP-semigroup $\tau$ on $\cB$, (the proof of) Theorem \ref{sgSkethm} provides us with isometries $v_t\colon E\odot E_t\rightarrow E$ such that $T_t=v_t(\bullet\odot\xi_t)$, where $(E^\odot,\xi^\odot)$ is the GNS-constructions for $\tau$. Of course, if these $v_t$ are adjointable, we are done by establishing $T$ as a strict CPH-semigroup. An excellent way of making sure that the $v_t$ have adjoints, would be if we could show that they are actually unitaries. In that case, $T$ would even be a strict CPH$_0$-semigroup.

We leave apart the question of adjointability, when the $v_t$ are non surjective. (Anyway, the situation that the GNS-system of $\tau$ sits adjointably in the product system of some CPH-extension $\sT$ as in Theorem \ref{sgthm} is not very likely. But for full $E$ it would be a necessary condition. And, anyway, except that Theorem \ref{sgthm} does not give a criterion by ``looking alone at $T$'', together with Corollary \ref{sgcor} it gives already a quite comprehensive answer to most questions.) $v_t$ being surjective, means
\beq{\label{Tmin}
E
~=~
\cls\BCB{\,T_{t_1}(T_{t_2}(\ldots T_{t_n}(x)b_{n-1}\ldots)b_1)b_0\colon n\in\N,t_1+\ldots+t_n=t,b_i\in\cB,x\in E\,}.
}\eeq
Since $v_s(E\odot E_s)\supset v_s(v_t(E\odot E_t)\odot E_s)=v_{s+t}(E\odot E_{s+t})$ for whatever CP-H-extendable semigroup $T$, the right-hand side decreases with $t$. So, it is sufficient to require that \eqref{Tmin} holds for some $t_0>0$.

\bdefi \label{minCPHdefi}
A CP-H-extendable semigroup $T$ on a Hilbert \nbd{\cB}module $E$ ($E$ full or not, $\cB$ unital or not) is \hl{minimal}, if $T$ fulfills \eqref{Tmin} for some $t_0>0$.
\edefi

Note that if $E$ is full (so that $\tau$ is unique) and if $T$ is minimal, then also the GNS-system of $\tau$ ($\cB$ unital or not; see Remark \ref{non1strirem}) is necessarily full. We are now ready to characterize minimal CP-H-extendable semigroups (which, therefore, are also CPH$_0$-semigroups) on full Hilbert modules over unital \nbd{C^*}algebras.

\bthm\label{minCPH0thm}
Let $\tau$ be a CP-semigroup on a unital \nbd{C^*}algebra $\cB$, and denote by $(E^\odot,\xi^\odot)$ its GNS-system and cyclic unit. Let $E$ be a full Hilbert \nbd{\cB}module. Then the formula $T_t=v_t(\bullet\odot\xi_t)$ establishes a one-to-one correspondence between:
\begin{enumerate}
\item
Left dilations $v_t\colon E\odot E_t\rightarrow E$ of $E^\odot$ to $E$.

\item
Minimal CP-H-extendable semigroups $T$ on $E$ associated with $\tau$.
\end{enumerate}
In either case, $\vt$ with $\vt_t=v_t(\bullet\odot\id_t)v_t^*$ is the unique strict \nbd{E_0}semigroup on $\sB^a(E)$ making $\sT=\rtMatrix{\tau&T^*\\T&\vt}$ a CPH$_0$-extension of $T$.
\ethm

\proof
Let $v_t$ be a left dilation. Then the \nbd{\tau_t}maps $T_t:=v_t(\bullet\odot\xi_t)$ define a CPH$_0$-semigroup $T$ on $E$. Since
\beqn{
E_t
~=~
\cls\bCB{b_n\xi_{t_n}\ldots b_1\xi_{t_1}b_0\colon n\in\N,t_1+\ldots+t_n=t,b_i\in\cB},
}\eeqn
we see that
\beqn{
T_{t_1}(T_{t_2}(\ldots T_{t_n}(x)b_{n-1}\ldots)b_1)b_0
~=~
x\xi_{t_n}\ldots b_1\xi_{t_1}b_0
}\eeqn
is indeed total in $v_t(E\odot E_t)=E$. Conversely, if $T$ is CP-H-extendable semigroup, we know see that
\beqn{
v_t
\colon
x\odot\xi_{t_n}\ldots b_1\xi_{t_1}b_0
~\longmapsto~
T_{t_1}(T_{t_2}(\ldots T_{t_n}(x)b_{n-1}\ldots)b_1)b_0
}\eeqn
defines isometries fulfilling \eqref{assoc}, which are unitary if and only if $T$ is minimal. Of course, $v_t(x\odot\xi_t)=T_t(x)$ so that the two directions are inverses of each other. This shows the one-to-one correspondence.

Finally, if $\vt_t$ is another endomorphism of $\sB^a(E)$, making $\sT_t$ a CPH-extension of $T_t$, then by the argument preceding Corollary \ref{CP(-)Hlem}, we have $\vt_t(a)T_t(x)=T_t(ax)$. So,
\bmun{
\vt_t(a)T_{t_1}(T_{t_2}(\ldots T_{t_n}(x)b_{n-1}\ldots)b_1)b_0
~=~
T_{t_1}(\vt_{t-t_1}(a)T_{t_2}(\ldots T_{t_n}(x)b_{n-1}\ldots)b_1)b_0
\\
~=~
\ldots
~=~
T_{t_1}(T_{t_2}(\ldots T_{t_n}(ax)b_{n-1}\ldots)b_1)b_0,
}\emun
that is, $\vt_t=v_t(\bullet\odot\id_t)v_t^*$.\qed

\bob
Recall that left dilations of $E^\odot$ to $E$ give rise to strict \nbd{E_0}semigroups on $\sB^a(E)$ that are all in the same cocycle equivalence class, and that every element in that cocycle equivalence class arises from such a left dilation. But different left dilations may have the same \nbd{E_0}semigroup; see, again, \cite[Proposition 6.3]{Ske08p1}. But it is the left dilations that are in one-to-one correspondence with the minimal CP-H-extendable semigroups. This underlines, once more, the importance of the concept of left dilations of a product system in addition to that of \nbd{E_0}semigroups associated with that product system.
\eob

There is cocycle version of the uniqueness result for the construction in \cite{Ske12a} proved by using the left dilations in Theorem \ref{minCPH0thm}. We state it without proof.

\bcor\label{minCP-Hunicor}
Let $T$ and $T'$ be two minimal CP-H-extendable semigroups on the same (necessarily full) Hilbert module $E$ over the unital \nbd{C^*}algebra $\cB$.

Then $T$ and $T'$ are associated with the same CP-semigroup $\tau$ on $\cB$ if and only if there is a unitary left cocycle for $\vt$ satisfying $u_t\colon T_t(x)\mapsto T'_t(x)$.

Moreover, if $u_t$ exists, then it is determined uniquely and $\vt'_t=u_t\vt_t(\bullet)u_t^*$.
\ecor

\noindent
So, minimal CP-H-extendable semigroups on the same $E$ associated with the same $\tau$ are no longer unitarily equivalent, but \it{cocycle equivalent}. We leave apart the question, when two minimal CP-H-extendable on the same $E$ but to possibly different $\tau$ have cocycle equivalent $\vt$, that is, their $\tau$ have isomorphic GNS-systems. The equivalence induced among Markov semigroups by their GNS-systems has been examined in Bhat and Skeide \cite[Section 7]{BhSk00}. It leads to a diffent sort of cocycles.

\bob
CP-H-extendable semigroups associated with the same fixed CP-semigroup $\tau$ may be added (direct sum), and the sum of minimal ones is again minimal. So, even if $E$ is full, minimality does not fix $E$ and $T$ up to cocycle equivalence. There is no \it{a priori} reason why two different $E$ should be isomorphic.
\eob

\bex
Let $\tau=\id_\cB$ be the trivial semigroup. So, CP-H-extendable semigroups associated with $\tau$ are just the semigroups of (\it{a priori} not necessarily adjointable) isometries on Hilbert \nbd{\cB}modules. It follows that $T_t(x)b=T_t(xb)$ so that minimality means $T_t(E)=E$. In other words, minimal CP-H-extendable semigroups associated with $\id_\cB$ are precisely the unitary semigroups. Of course, (for suitable $\cB$, for instance, for $\cB=\C$) there are nonisomorphic full Hilbert \nbd{\cB}modules.
\eex

However, if $E$ is full and countably generated (over unital $\cB$, so that $\cB$ is in particular \nbd{\sigma}unital), then $E^\infty\cong\cB^\infty$; see Lance \cite[Proposition 7.4]{Lan95}. So, minimal CP-H-extendable semigroups on different countably generated full $E$ may, first, be lifted to $\cB^\infty$ and, then, there are cocycle equivalent. In other words, the original semigroups are \it{stably} cocycle equivalent.

\bob \label{notETob}
Whatever the CP-H-extendable semigroup $T$ is, if $\cB$ is unital, then
\beqn{
E_{T_t}
~=~
\cls T_t(E)\cB
~=~
\cls v_t(E\odot\xi_t\cB)
~=~
v_t(E\odot\sF_t),
}\eeqn
where $\sF_t=\cls\cB\xi_t\cB\subset E_t$ is the GNS-correspondences of the single CP-map $\tau_t$. It is a typical feature of the GNS-system that $\sF_{s+t}\subset\cls\sF_s\sF_t\subset E_{s+t}$ is smaller than $E_{s+t}$ unless $\tau$ is an \nbd{E}semigroup, because $\sF_{s+t}\cong\cls\cB\xi_s\odot\xi_t\cB\subsetneq\cls\cB\xi_s\odot\cB\xi_t\cB=\sF_s\odot\sF_t$.
\eob

\brem
If $T$ is not minimal, then the ranges of the $v_t$ decrease, say, to $E_\infty$. Moreover, it is clear that the $v_t$ (co)restrict to unitaries $E_\infty\odot E_t\rightarrow E_\infty$, that is, they form a left quasi-dilation of $E^\odot$ to $E_\infty$. It is unclear if $E_\infty$ is full, even if $E$ is full and $E^\odot$ is full, or if $E_\infty$ may be possibly $\zero$. But in any case, $T$ (co)restricts to a minimal strictly CPH$_0$-semigroup on $E_\infty$ associated with $\tau$. Necessarily, $\tau$ (co)restricts to a CP-semigroup on $\cB_{E_\infty}$.
\erem

\lf
It might be worth to compare the results in this section with Heo and Ji \cite{HeJi11}, who investigated semigroups that, in our terminology, are CP-H-extendable, but who call them CP-semigroups.


\lf\lf
\section{An application: CPH-dilations}\label{CPHdilSEC}


Since Asadi drew attention to \nbd{\tau}maps $T\colon E\rightarrow\sB(G,H)$ for CP-maps $\tau\colon\cB\rightarrow\sB(G)$, it is an open question what they might be good for. In this section, we make the first attempt to give them an interpretation; and our point is to interpret them as a notion that generalizes the notion of dilation of a CP-map $\tau\colon\cB\rightarrow\cC$ to a homomorphism $\vt\colon\sB^a(E)\rightarrow\sB^a(F)$ to the notion of \it{CPH-dilation}. In particular, in the situation of semigroups, our new more relaxed version of dilation allows for new features: While CP-semigroups that allow \it{weak} dilations to an \nbd{E_0}semigroup (also \it{\nbd{E_0}dilations}), are necessarily Markov, our results from Section \ref{CPHsgSEC} allow us to show that many nonunital CP-semigroups allow CPH-dilations to \nbd{E_0}semigroups, \it{CPH$_0$-dilations}.

Let us start with a CP-map $\tau\colon\cB\rightarrow\cC$ with unital $\cB$, and with a \nbd{\tau}map $T\colon E\rightarrow F$. Denoting by $(\sF,\zeta)$ the GNS-construction for $\tau$, by \cite{Ske12a} we get a (unique) isometry $v\colon E\odot\sF\rightarrow F$ such that $T(x)=v(x\odot\zeta)$. If $F_T$ is complemented in $F$, that is, if $v$ is adjointable, then $\vt\colon a\mapsto v(a\odot\id_\sF)v^*$ is a strict homomorphism from $\sB^a(E)$ to $\sB^a(F)$. Now, if $\xi$ is a unit vector (that is, $\AB{\xi,\xi}=\U$) in $E$, we may define the representation $b\mapsto\xi b\xi^*$ of $\cB$ on $E$. We find
\beq{\label{dileq}
\AB{v(\xi\odot\zeta),\vt(\xi b\xi^*)v(\xi\odot\zeta)}
~=~
\AB{\xi\odot\zeta,(\xi b\xi^*\odot\id_\sF)(\xi\odot\zeta)}
~=~
\AB{\zeta,b\zeta}
~=~
\tau(b),
}\eeq
so that the following diagram commutes.
\beqn{
\xymatrix{
\cB	 		\ar[rr]^{\tau}		\ar[d]_{\xi\bullet\xi^*}	&&	\cC
\\
\sB^a(E)	\ar[rr]_{\vt}											&&	\sB^a(F)		\ar[u]_{\AB{v(\xi\odot\zeta),\bullet v(\xi\odot\zeta)}}
}
}\eeqn
It is clear that just any quintuple $(\sF,\zeta,E,v,\xi)$ of a \nbd{\cB}\nbd{\cC}correspondence $\sF$, an element $\zeta\in\sF$, a Hilbert \nbd{\cB}module $E$, an adjointable isometry $v\colon E\odot\sF\rightarrow F$, and a unit vector $\xi\in E$ will do, if we put $\tau:=\AB{\zeta,\bullet\zeta}$ and $\vt:=v(\bullet\odot\id_\sF)v^*$.

If also $\zeta$ is a unit vector (so that $\tau$ is unital, and also $v(\xi\odot\zeta)$ is a unit vector), such a situation is called a \hl{weak dilation} of the \hl{Markov map} (that is, a unital CP-map) $\tau$. Here `weak' is referring to that the embedding $\cB\rightarrow\xi\cB\xi^*$ means identifying $\cB$ with a corner in $\sB^a(E)$ (and likewise $\cC\rightarrow(\xi\odot\zeta)\cC(\xi\odot\zeta)^*$) and that $\xi\xi^*\bullet\xi\xi^*=\xi\AB{\xi,\bullet\xi}\xi^*$ is just the conditional expectation onto that corner (and likewise for the corner of $\sB^a(F)$ isomorphic to $\cC$).

What, if we do not have a unit vector in $E$ or if $\tau$ is not unital? Let us make two observations: Firstly, as long as $\xi$ is a unit vector, the condition that the preceding diagram commutes is actually equivalent to the apparently stronger condition that the diagram
\beqn{
\xymatrix{
\cB	 		\ar[rr]^{\tau}												&&	\cC
\\
\sB^a(E)	\ar[rr]_{\vt}		\ar[u]^{\AB{\xi,\bullet\xi}}	&&	\sB^a(F)	\ar[u]_{\AB{v(\xi\odot\zeta),\bullet v(\xi\odot\zeta)}}	
}
}\eeqn
commutes. For this, $\zeta$ need not be a unit vector. (In fact, substituting in \eqref{dileq} $\xi b\xi^*$ with $a\in\sB^a(E)$, the same computation yields $\AB{v(\xi\odot\zeta),\vt(a)v(\xi\odot\zeta)}=\tau(\AB{\xi,a\xi})$, and inserting $a=\xi b\xi^*$ gives back the original equation.) Secondly, in the expectation the \nbd{\tau}map $T:=v(\bullet\odot\zeta)$ occurs as $\AB{v(\xi\odot\zeta),\bullet v(\xi\odot\zeta)}=\AB{T(\xi),\bullet T(\xi)}$. In this form, the diagram makes sense also if we replace the $\xi$ in the left factor and the $\xi$ in the right factor of the inner products with an arbitrary pair $x,x'$ of elements of $E$:

\bdefi\label{CPHdildefi}
Let $\tau\colon\cB\rightarrow\cC$ be a CP-map. A homomorphism $\vt\colon\sB^a(E)\rightarrow\sB^a(F)$ is a \hl{CPH-dilation} of $\tau$ if $E$ is full and if there is a map $T\colon E\rightarrow F$ such that the diagram
\beqn{
\xymatrix{
\cB	 		\ar[rr]^{\tau}												&&	\cC
\\
\sB^a(E)	\ar[rr]_{\vt}		\ar[u]^{\AB{x,\bullet x'}}	&&	\sB^a(F)		\ar[u]_{\AB{T(x),\bullet T(x')}}
}
}\eeqn
commutes for all $x,x'\in E$. (We do not require that $\cB$ and $\cC$ are unital.) If $E$ is not necessarily full, then we speak of a \hl{CPH-quasi-dilation}. A CPH-(quasi)dilation is \hl{strict} if $\vt$ is strict. A CPH-(quasi-)dilation is a \hl{CPH$_0$-}(\hl{quasi-})\hl{dilation} if $\vt$ is unital.
\edefi

Requiring dilation instead of quasi-dilation, means excluding trivialities. (Without that, $E$ may be very well $\zero$.) Of course, a CPH-dilation may be turned into a CPH$_0$-dilation, by replacing $F$ with $\vt(\id_E)F$. It is strict if and only if  $\vt(\id_E)F=\cls\vt(\sK(E))F$. In a CPH-quasi-dilation, the diagram does not give any information about the component of $T(x)$ in $(\vt(\id_E)F)^\perp$. In that case, it is convenient to replace $T$ with $\vt(\id_E)T$ and apply the following results to the latter map considered as map into $\vt(\id_E)F$.

\bprop\label{llinTprop}
If $\vt$ is a CPH$_0$-quasidilation of a CP-map $\tau$, then every map $T$ making the diagram commute is a \nbd{\tau}map fulfilling $T(ax)=\vt(a)T(x)$.
\eprop

\proof
Inserting $a=\id_E$ into the diagram, we see that $T$ is a \nbd{\tau}map. Also, for arbitrary $a\in\sB^a(E)$ and $x,x'\in E$, we get $\AB{T(x),\vt(a)T(x')}=\AB{T(x),T(ax')}$. A brief argument shows that this implies $\vt(a)T(x)=T(ax)$. (Indeed, on $F_T:=\cls T(E)\cC\subset F$, we know we get a (strict, unital) representation $\vt_T\colon\sB^a(E)\rightarrow\sB^a(F_T)$ that acts on the generating subset $T(E)$ in the stated way. That is, we have $\AB{y,\vt(a)y'}=\AB{y,\vt_T(a)y'}$ for all $y,y'\in F_T$. From this, one easily verifies that $\abs{\vt(a)y-\vt_T(a)y}^2=0$, so, $\vt(a)y=\vt_T(a)y$ for all $y\in F_T$.)\qed

\lf
The \nbd{\vt}left linearity of $T$ looks like something we would knew already from Lemma \ref{bilcp} and the discussion following it. Note, however, that this discussion is based entirely on the assumption that the extension $\sT=\rtMatrix{\tau&T^*\\T&\vt}$ is a CP-map --- a hypothesis we still do not yet know to be true. In fact, we will prove it in the following theorem only for strict CPH$_0$-dilations for unital $\cB$. And still there it turns out to be surprisingly tricky.

\lf
From now on we shall assume that $\cB$ is \bf{unital}.

\bthm\label{CPHdilthm}
If $\vt$ is a strict CPH$_0$-dilation of a CP-map $\tau$, then every map $T$ making the diagram commute is a strict CPH$_0$-map.
\ethm

\proof
We shall show that $\sT=\rtMatrix{\tau&T^*\\T&\vt}$ is CP, so that $T$ is strictly CPH$_0$. We wish to imitate the proof of complete positivity in the step \ref{1} $\Rightarrow$ \ref{2} in Section \ref{proofSEC}. But we have to face the problem that the multiplicity correspondence of $\vt$ does no longer coincide with the GNS-correspondence of $\tau$; it just contains it.

Denote by $(\sF_\tau,\zeta)$ the GNS-construction for $\tau$. Doing the representation theory for the strict unital homomorphism $\vt$, we get a \nbd{\cB}\nbd{\cC}correspondence $\sF_\vt:=E^*\odot{_\vt}F$ and a unitary $v\colon E\odot\sF_\vt\rightarrow F, x'\odot(x^*\odot y)\mapsto\vt(x'x^*)y$ such that $\vt=v(\bullet\odot\id_{\sF_\vt})v^*$. By Proposition \ref{llinTprop}, one easily verifies that
\beqn{
\AB{x,x'}\zeta
~\longmapsto~
x^*\odot T(x')
}\eeqn
determines a bilinear unitary from $\sF_\tau$ onto $E^*\odot F_T\subset\sF_\vt$. We shall identify $\sF_\tau\subset\sF_\vt$, so that $\zeta\in\sF_\vt$. We have
\beqn{
v(x\odot\AB{x',x''}\zeta)
~=~
v(x\odot(x'^*\odot T(x''))
~=~
\vt(xx'^*)T(x'')
~=~
T(x\AB{x',x''}).
}\eeqn
Since $\ls\AB{E,E}\ni\U$ and $v$ and $T$ are linear, it follows that $T(x)=v(x\odot\zeta)$.

Since $F_T$ need not be complemented in $F$, also $\sF_\tau$ need not be complemented in $\sF_\vt$. But, we still have a map $\rtMatrix{\zeta~&\\&~v^*}\in\sB^r\rtMatrix{\tMatrix{\cC\\F},\tMatrix{\cB\\E}\odot\sF_\vt}$. We find
\beq{ \label{v-T}
\SMatrix{\SMatrix{b&x^*\\x'&a}\odot\id_{\sF_\vt}}\SMatrix{\zeta~&\\&~v^*}\SMatrix{c\\y}
~=~
\SMatrix{b\zeta c+(x^*\odot\id_{\sF_\vt})v^*y\\x'\odot\zeta c+(a\odot\id_{\sF_\vt})v^*y},
}\eeq
so that
\baln{
\BbAB{
\SMatrix{\SMatrix{b_1&x_1^*\\x'_1&a_1}\odot\id_{\sF_\vt}}\SMatrix{\zeta~&\\&~v^*}\SMatrix{c_1\\y_1},
\SMatrix{\SMatrix{b_2&x_2^*\\x'_2&a_2}\odot\id_{\sF_\vt}}\SMatrix{\zeta~&\\&~v^*}\SMatrix{c_2\\y_2}
}\hspace{-30ex}
\\
~=~
&
c_1^*\AB{\zeta,b_1^*b_2\zeta}c_2
+c_1^*\AB{\zeta,b_1^*(x_2^*\odot\id_{\sF_\vt})v^*y_2}
\\
&~~~~~~~~~~~~
+\AB{(x_1^*\odot\id_{\sF_\vt})v^*y_1,b_2\zeta}c_2
+\AB{(x_1^*\odot\id_{\sF_\vt})v^*y_1,(x_2^*\odot\id_{\sF_\vt})v^*y_2}
\\
&
+c_1^*\AB{x'_1\odot\zeta,x'_2\odot\zeta}c_2
+c_1^*\AB{x'_1\odot\zeta,(a_2\odot\id_{\sF_\vt})v^*y_2}
\\
&~~~~~~~~~~~~
+\AB{(a_1\odot\id_{\sF_\vt})v^*y_1,x'_2\odot\zeta}c_2
+\AB{(a_1\odot\id_{\sF_\vt})v^*y_1,(a_2\odot\id_{\sF_\vt})v^*y_2}
\\
~=~
&
c_1^*\tau(b_1^*b_2)c_2
+c_1^*\AB{T(x_2b_1),y_2}
+\AB{y_1,T(x_1b_2)}c_2
+\AB{y_1,\vt(x_1x_2^*)y_2}
\\
&
+c_1^*\tau(\AB{x'_1,x'_2})c_2
+c_1^*\AB{T(a_2^*x'_1),y_2}
+\AB{y_1,T(a_1^*x'_2)}c_2
+\AB{y_1,\vt(a_1^*a_2)y_2}
\\
&~~~~~~~~~~~~~~~~~~~~~~~~~~~~~~~~~~~~~~~~~~~~~~~~~~~~~~~~~~~~
~=~
\BbAB{\SMatrix{c_1\\y_1},\sT\SMatrix{\SMatrix{b_1&x_1^*\\x'_1&a_1}^*\SMatrix{b_2&x_2^*\\x'_2&a_2}}\SMatrix{c_2\\y_2}}.
}\ealn
Taking appropriate sums of such expressions, we see that $\sT$ is completely positive.\qed

\bob
It is crucial that we define an embedding from $\sF_\tau$ into $\sF_\vt$ by fixing its values on $\AB{x,x'}\zeta$. Only if $E$ is full, this determines an isometry on all of $\sF_\tau$. And to be sure $\zeta$ exists, $\cB$ has to be unital.

If $\cB$ is nonunital (still $E$ full), then instead of $\zeta$ we may look at elements $\AB{x,x'}\otimes c+\sN$ in the GNS-correspondence $\sF_\tau=\ol{(\cB\otimes\cC)/\sN}$. We define $\sF_\tau\rightarrow\sF_\vt$ as $\AB{x,x'}\otimes c+\sN\mapsto x^*\odot T(x')c$. Instead of \eqref{v-T}, we consider the elements
\beqn{
\lim_\lambda\SMatrix{(bu_\lambda\otimes c+\sN)+(x^*\odot\id_{\sF_\vt})v^*y\\x'\odot(u_\lambda\otimes c+\sN)+(a\odot\id_{\sF_\vt})v^*y}
}\eeqn
in $\tMatrix{\cB\\E}\odot\sF_\vt$, where $\bfam{u_\lambda}_{\lambda\in\Lambda}$ is an approximate unit in $\ls\AB{E,E}$ for $\cB$. Everything in the long computation of the proof of Theorem \ref{CPHdilthm} goes through as before, showing that $T$ is a strictly CPH$_0$-map. But in this experimental section we do not intend to be exhaustive, and stick to the simplest case where $E$ is full over unital $\cB$.
\eob

Appealing to Theorem \ref{mainthm}, \ref{3} $\Rightarrow$ \ref{2}, and Observation \ref{discob}\eqref{O4}, we get following:

\bcor
If $E$ is full over unital $\cB$ and $\vt\colon\sB^a(E)\rightarrow\sB^a(F)$ is a strict unital homomorphism for which there exists a linear map $T\colon E\rightarrow F$ such that $T(ax)=\vt(a)T(x)$, then each such $T$ is a strict CPH$_0$-map and $\vt$ is a strict CPH$_0$-dilation of the CP-map $T^*\odot T$.
\ecor

Note that, without fixing $\tau$, every homomorphisms $\vt$ is a CPH-dilation of the CP-map $\tau=0$. So, CPH-dilation is meaningful only with reference to a fixed CP-map.

$T$ need not be unique, not even up to unitary equivalence.

\bex\label{nonuniex}
Let $F$ be such that $F\oplus F\cong F$ as \nbd{\sB^a(E)}\nbd{\cC}correspondences. If $T$ is good enough to make the diagram commute, then so is either map $T_i$ sending $x$ to $T(x)$ in the \nbd{i}component of $F\oplus F$. Essentially, given $\tau$ and $\vt$, it is undetermined how $F_T$ sits inside $F$ and even $F_T^\perp$ for different $T$ need not be isomorphic.
\eex

However, as usual, if we require, for given $F$, that the map $T$ fulfills $F_T=F$, then we know that up to unitary automorphism $u$ of $F$ leaving $\vt$ invariant, there is at most one $T$. Note that the unitaries on $F_T\cong E\odot\sF$ not changing $\vt$, have to commute with $\vt(\sB^a(E))\cong\sB^a(E)\odot\id_\sF$. For full $E$, this means $u\cong\id_E\odot\upsilon$ for some automorphism $\upsilon\in\sB^{a,bil}(\sF)$ of the GNS-correspondence $\sF$ of $\tau$.

We see, different minimal $T$ are distinguished by ``shoving around'' (with $\upsilon$) the cyclic vector $\zeta$ that occurs in $T(x)=v(x\odot\zeta)$. But doing so, under minimality, we get unitarily equivalent things. This gets much more interesting in the semigroup case, to which we switch now, where this ``shoving around'' has to be done compatibly with the semigroup structure. Recall that even for a single CP-map $\tau$ not \bf{between} \nbd{C^*}algebras but \bf{on} a \nbd{C^*}algebra, it is required that the dilating map $\vt$ does not only dilate $\tau$, but that for each $n$ the power $\vt^{\circ n}$ dilates the power $\tau^{\circ n}$. In particular, we will see that the usual concept of \it{weak dilation} of a CP-semigroup (of which CPH-dilations are a generalization) means that the corresponding semigroup $T$ has to leave the vector $\xi$ fixed.

Let us begin with this situation of weak dilation, by continuing the report on results from Bhat and Skeide \cite{BhSk00}. We mentioned already in Section \ref{CPHsgSEC} that for every CP-semigroup $\tau$ on a unital \nbd{C^*}algebra $\cB$, we get the GNS-construction $(E^\odot,\xi^\odot)$ consisting of a product system $E^\odot$ and unit $\xi^\odot$ for $E^\odot$ that generates $E^\odot$ and that gives back $\tau$ as $\tau_t=\AB{\xi_t,\bullet\xi_t}$. The semigroup $\tau$ is Markov if and only if the unit $\xi^\odot$ is \bf{unital}, that is, if $\AB{\xi_t,\xi_t}=\U$ for all $t$. Starting from a product system with a unital unit, \cite{BhSk00} provide the following additional ingredients:
\begin{itemize}
\item
A left dilation $v_t\colon E\odot E_t\rightarrow E$ of $E^\odot$ to a (by definition full) Hilbert module $E$. So, the maps $\vt_t\colon a\mapsto v_t(a\odot\id_t)v_t^*$ define a strict \nbd{E_0}semigroup on $\sB^a(E)$.

\item
A unit vector $\xi\in E$ such that $\xi\xi_t=\xi$. It is readily verified that the triple $(E,\vt,\xi)$ is a \hl{weak dilation} of $\tau$ in the sense that
\beqn{
\AB{\xi,\vt_t(\xi b \xi^*)\xi}
~=~
\tau_t(b).
}\eeqn
In other words, if we define the projection $p:=\xi\xi^*\in\sB^a(E)$ and identify $\cB$ with the corner $\xi\cB\xi^*=p\sB^a(E)p$ of $\sB^a(E)$, then $p\vt_t(a)p=\tau_t(pap)\in\sB^a(E)$.

If $(E^\odot,\xi^\odot)$ is the GNS-construction, then the dilation constructed in\cite{BhSk00} is \hl{minimal} in the sense that $\vt_{\R_+}(\xi\cB\xi^*)$ generates $E$ out of $\xi$. Such a minimal dilation is unique up to suitable unitary equivalence.
\end{itemize}
Now, if we define $T_t(x):=x\xi_t$, we see that the diagram
\beq{\label{CPHSGdiag}
\begin{minipage}{6cm}
\xymatrix{
\cB	 		\ar[rr]^{\tau_t}												&&	\cB
\\
\sB^a(E)	\ar[rr]_{\vt_t}		\ar[u]^{\AB{x,\bullet x'}}	&&	\sB^a(E)		\ar[u]_{\AB{T_t(x),\bullet T_t(x')}}
}
\end{minipage}
}\eeq
commutes for all $x,x'\in E$ and all $t\in\R_+$. The special property of the dilation from \cite{BhSk00} is the existence of the unit vector $\xi\in E$ fulfilling $\xi\xi_t=\xi$, that is, $T_t$ leaves $\xi$ fixed. But for that the diagram commutes, effectively just any left dilation $v_t$ will do. For any product system $E^\odot$, any unit $\xi^\odot$ and any left dilation $v_t\colon E\odot E_t\rightarrow E$ to a full Hilbert \nbd{\cB}module $E$ (so that all $E_t$ are necessarily full, too), the formulae $\tau_t:=\AB{\xi_t,\bullet\xi_t}$, $\vt_t:=v_t(\bullet\odot\id_t)v_t^*$, and $T_t(x):=x\xi_t$ provide us with a strict CPH$_0$-dilation of $\tau_t$. For this it is not necessary that the $\tau_t$ form a Markov semigroup. Of course, also the corresponding $\sT_t$ form a (strict) CP-semigroup (which is Markov if and only if $\tau_t$ is Markov).

\bdefi
An \nbd{E_0}semigroup $\vt$ on $\sB^a(E)$ for a full Hilbert \nbd{\cB}module $E$ is a \hl{CPH$_0$-dilation} of a CP-semigroup $\tau$ on $\cB$ if there exists a CPH$_0$-semigroup $T$ on $E$ making Diagram \eqref{CPHSGdiag} commute for all $t\in\R_+$. (We use all variants as in Definition \ref{CPHdildefi}.)
\edefi

If $\tau_t$ is not Markov, then \cite{BhSk00} provide a weak dilation to an \nbd{E}semigroup. But $\tau_t$ cannot posses a weak dilation to an \nbd{E_0}semigroup. On the contrary, we see that $\tau_t$ \bf{can} possess a CPH$_0$-dilation:

\bob\label{CPHdilob}
Finding a strict CPH($_0$)-dilation for a CP-semigroup $\tau$, is the same as finding a CPH($_0$)-semigroup $T$ associated with that $\tau$. So, all our results from Section \ref{CPHsgSEC} are applicable.

\begin{enumerate}
\item
From Corollary \ref{sgcor}, we recover existence of a strict CPH-dilation. (As said, we knew this from the stronger existence of a weak dilation in \cite{BhSk00}.)

\item
But, in particular, as in the discussion following Corollary \ref{sgcor}, from existence of \nbd{E_0}semi\-groups for full product systems, we infer that every CP-semigroup, Markov or not, with full product system admits a strict CPH$_0$-dilation.

\item
In the case of CPH$_0$-dilations, also the notion of minimality and the results about uniqueness up to cocycle conjugacy remain intact. It is noteworthy that for a weak \nbd{E_0}dilation of a (necessarily) Markov semigroup, minimality of the weak dilation coincides with minimality of the associated CPH$_0$-semigroup.
\end{enumerate}
\eob

\lf\noindent
We see that CPH-semigroups and CPH-dilations are to some extent two sides of the same coin --- a coin that can be expressed as in the diagram of CPH-dilation in \eqref{CPHSGdiag}. CPH-maps put emphasis of the map between the modules, and under suitable cyclicity requirements the remaining corners $\tau$ (if $E$ is full) and $\vt$ (if $F_T=F$) follow. CPH-dilations put emphasis on that there is a \it{relation} between the diagonal corners. While the notion of CPH-dilation underlines that we are in front of a generalized dilation of a CP-semigroup to an endomorphism semigroup (namely, where there is no longer a cyclic vector, and if it is there it need no longer be fixed by the associated CPH-semigroup), the notion of CPH-semigroup underlines that there is, at least under good cyclicity conditions, a single object, the CPH-semigroup, that encodes everything and that may be studied separately.

We close by some considerations regarding situations related with CPH-di\-la\-tions, which might be interesting. This is not any concrete evidence, but for now mere speculation. But if some of these situations, in the future, really will turn out to be interesting, the mutual relation between CPH-dilations and CPH-semigroups, in particular the results of Section \ref{CPHsgSEC}, will find their applications. After all, while so far all publications about CPH-maps and CP-H-extendable semigroups are justified only by claiming interest ``on their own'', our considerations here, though rather speculative, are the first pointing into the direction of potential applications.

\lf
Let us have a different look at Diagram \eqref{CPHSGdiag}. Note that the map $\eK\colon(x,x')\mapsto\eK^{x,x'}:=\AB{x,\bullet x'}$ is a \hl{completely positive definite} or \hl{CPD-kernel} over the set $E$ from $\sB^a(E)$ to $\cB$ in the sense of Barreto, Bhat, Liebscher, and Skeide \cite[Section 3.2]{BBLS04}; see also the survey Skeide \cite{Ske11}. The maps $T_t$ amount to a transformation semigroup of the indexing set $E$. We may generalize CPH-dilation of a CP-semigroup $\tau$ on $\cB$ to the situation
\beqn{
\xymatrix{
\cB	 		\ar[rr]^{\tau_t}												&&	\cB
\\
\cA	\ar[rr]_{\theta_t}		\ar[u]^{\eK^{\sigma,\sigma'}}	&&	\cA		\ar[u]_{\eK^{T_t(\sigma),T_t(\sigma')}}
}
}\eeqn
where $\theta$ is an endomorphism semigroup on a unital \nbd{C^*}algebra $\cA$ and where $\eK$ is a fixed CPD-kernel over $S$ from $\cA$ to $\cB$. Note, however, that this situation is not too much more general. Effectively, $\eK$ has a Kolmogorov decomposition $(E,\vk)$ consisting of an \nbd{\cA}\nbd{\cB}correspondence $E$ and a map $\vk\colon S\rightarrow E$ such that $\eK^{\sigma,\sigma'}=\AB{\vk(\sigma),\bullet\vk(\sigma')}$ and $E=\cls\cA\vk(S)\cB$.

A natural question is if $T_t$ extends as a map $E\rightarrow E$ (automatically a \nbd{\tau_t}map). Another question is if $\cA$ is $\sB^a(E)$, and, if not, if there is an \nbd{E}semigroup $\vt$ on $\sB^a(E)$ such that the left action of $\theta_t(a)$ on $E$ is the same as $\vt_t$ applied to the operator on $E$ given by the left action of $a$. (These questions are direct generalizations of the same questions for usual dilations of CP-semigroups: Does every dilation to $\cA$ give rise to a dilation to $\sB^a(E)$ where $E$ is the GNS-correspondence of the conditional expectation?)

We also may ask, if this setting has a useful interpretation in terms of Morita equivalence. If $\cA=\sB^a(E)$ and if $E$ is full, then $\sK(E)$ is Morita equivalent to $\cB$. We may say, $\sB^a(E)$ is strictly Morita equivalent to $M(\cB)$. The CPD-kernel somehow encodes the necessary information about the Morita equivalence transform: The identification $\cB=E^*\odot E=E^*\odot\sB^a(E)\odot E$, giving rise to the kernel $\eK^{x,x'}(a)=x^*\odot a\odot x'=\AB{x,ax'}$. How is the transform $T_t$ reflected in the picture of Morita equivalence? Is the Kolmogorov construction for $\eK^{T_t(\sigma),T_t(\sigma')}$ in a reasonable way contained in $E$? Of course, Morita equivalence is invertible. Is the ``inverse'' CPD-kernel $\eL$ from $\cB$ to $\cA$ defined by $\eL^{x',x}(b):=(x'^*)^*\odot b\odot x^*=x'bx^*$ of any use?

Answers to these questions will have to wait for future investigation.

\lf\lf\lf\noindent
\bf{Acknowledgments:~}
This work grew out of a six months stay during the first named author's sabbatical in 2010 and a three months stay in 2012 at ISI Bangalore. He wishes to express his gratitude to Professor Bhat and the ISI Bangalore for warm hospitality and the Dipartimento S.E.G.e S.\ of the University of Molise as well as the Italian MIUR (PRIN 2007) for taking over travel expenses. The second named author expresses his gratitude to Professor Bhat for his gracious encouragement.

We wish to thank Orr Shalit and Harsh Trivedi who contributed very useful comments on a preliminary version.\lf


\setlength{\baselineskip}{2.5ex}


\newcommand{\Swap}[2]{#2#1}\newcommand{\Sort}[1]{}
\providecommand{\bysame}{\leavevmode\hbox to3em{\hrulefill}\thinspace}
\providecommand{\MR}{\relax\ifhmode\unskip\space\fi MR }
\providecommand{\MRhref}[2]{%
  \href{http://www.ams.org/mathscinet-getitem?mr=#1}{#2}
}
\providecommand{\href}[2]{#2}

\lf\lf\lf\noindent
Michael Skeide,
{\small\itshape Universit\`a\ degli Studi del Molise},
{\small\itshape Dipartimento S.E.G.e S.},
{\small\itshape Via de Sanctis}
{\small\itshape 86100 Campobasso, Italy},
{\small{\itshape E-mail: \tt{skeide@unimol.it}}},\\
{\small{\itshape Homepage: \tt{http://web.unimol.it/skeide/}}}.

\lf\noindent
K.\ Sumesh,
{\small\itshape Indian Statistical Institute Bangalore},
{\small\itshape R.V.\ College Post},
{\small\itshape 8th Mile Mysore Road}
{\small\itshape 560059 Bangalore, India},
{\small{\itshape E-mail: \tt{sumesh@isibang.ac.in}}}.


\end{document}